\documentclass[times]{article}

\usepackage{moreverb}

\usepackage{tikz,pgfplots}
\usetikzlibrary{positioning,calc, shapes, backgrounds}

\usepackage{amsmath}
\usepackage{amssymb}
\usepackage{amsbsy,braket}
\usepackage[titletoc]{appendix}
\usepackage{booktabs}
\usepackage{mathrsfs}
\usepackage{todonotes,subfigure}
\usepackage{url}
\usepackage{fullpage}
\tikzset{>=latex}

\newcommand{\bd}{\partial}

\renewcommand{\gg}{\mathbf{g}}

\newcommand{\ndotnx}{\nn \cdot \nn_{\xx}}
\newcommand{\nn}{\mathbf{n}}
\newcommand{\ttau}{\boldsymbol{\tau}}
\newcommand{\pderiv}[2]{\frac{\partial #1}{\partial #2}}
\newcommand{\rdotn}{\rr \cdot \nn}
\newcommand{\rdotnx}{\rr \cdot \nn_{\xx}}
\newcommand{\rr}{\mathbf{r}}

\newcommand{\ssigma}{\boldsymbol{\sigma}}

\newcommand{\xx}{\mathbf{x}}
\newcommand{\yy}{\mathbf{y}}

\newcommand{\eps}{{\displaystyle \varepsilon}}
\newcommand{\bsub}{\begin{subequations}}
\newcommand{\esub}{\end{subequations}$\!$}

\newcommand{\bigoh}{\mathcal{O}}

\begin{document}

\title{A boundary integral equation method for mode elimination and
vibration confinement in thin plates with clamped points}
\author{Alan E. Lindsay, Bryan Quaife, and Laura Wendelberger}

\maketitle

\begin{abstract}
We consider the bi-Laplacian eigenvalue problem for the modes of
vibration of a thin elastic plate with a discrete set of clamped points.
A high-order boundary integral equation method is developed for
efficient numerical determination of these modes in the presence of
multiple localized defects for a wide range of two-dimensional
geometries. The defects result in eigenfunctions with a weak singularity
that is resolved by decomposing the solution as a superposition of
Green's functions plus a smooth regular part. This method is applied to
a variety of regular and irregular domains and two key phenomena are
observed. First, careful placement of clamping points can entirely
eliminate particular eigenvalues and suggests a strategy for
manipulating the vibrational characteristics of rigid bodies so that
undesirable frequencies are removed.  Second, clamping of the plate can
result in partitioning of the domain so that vibrational modes are
largely confined to certain spatial regions. This numerical method gives
a precision tool for tuning the vibrational characteristics of thin
elastic plates. 
\end{abstract}

%%%%%%%%%%%%%%%%%%%%%%%%%%%%%%%%%%%%%%%%%%%%%%%%%%%%%%%%%%%%%%%%%%%%%%%%
\section{Introduction}
\label{sec:introduction}
%%%%%%%%%%%%%%%%%%%%%%%%%%%%%%%%%%%%%%%%%%%%%%%%%%%%%%%%%%%%%%%%%%%%%%%%
The eigenvalues of fourth-order differential operators are central in
determining mechanical properties of rigid bodies. This paper considers
the small amplitude out-of-plane vibrations of a thin elastic
plate~\cite{RE}. The vibrational frequencies $\lambda>0$ and modes
$u(\xx)$ satisfy the bi-Laplacian eigenvalue problem
\bsub\label{eqn:intro}
  \begin{equation}
  \label{eqn:introA} \Delta^2 u = \lambda u, \qquad \xx \in\Omega; \qquad  \int_{\Omega} u^2\, d\xx = 1,
  \end{equation}
where $\Omega \subset \mathbb{R}^2$ is a closed planar region
representing the extent of the plate, $\xx = (x,y)$, and $\Delta^2 u:=
u_{xxxx} + 2 u_{xxyy} + u_{yyyy}$. Conditions on the boundary
$\partial\Omega$ are application specific, with a common condition
being that the plate is \emph{clamped} on its periphery which stipulates that 
  \begin{equation}\label{eqn:introB} 
  u = \partial_\nn u = 0, \qquad \xx \in\partial\Omega,
  \end{equation} 
where $\partial_{\nn}$ is the outward facing normal derivative. A wide
variety of engineering systems utilize thin perforated plates in their
construction. Examples include heat exchangers~\cite{Nilles95,
Venkatarathnam1996, Krishnakumar2003}, porous elastic materials, and
acoustic tilings~\cite{Atalla2007, wang2010, Jaouen2011}. The specific
placement of these perforations permits the manipulation of acoustic
and vibrational properties of the plate while economizing on weight and
material cost. Homogenization theories have been proposed to replace
the natural elastic modulus of the plate with an effective
modulus~\cite{BH,ADK}, however, an averaging approach omits the
pronounced localizing effects that clamping has on vibrational
modes~\cite{FM}.

\begin{figure}[htbp]
\centering
\includegraphics[width = 0.9\textwidth]{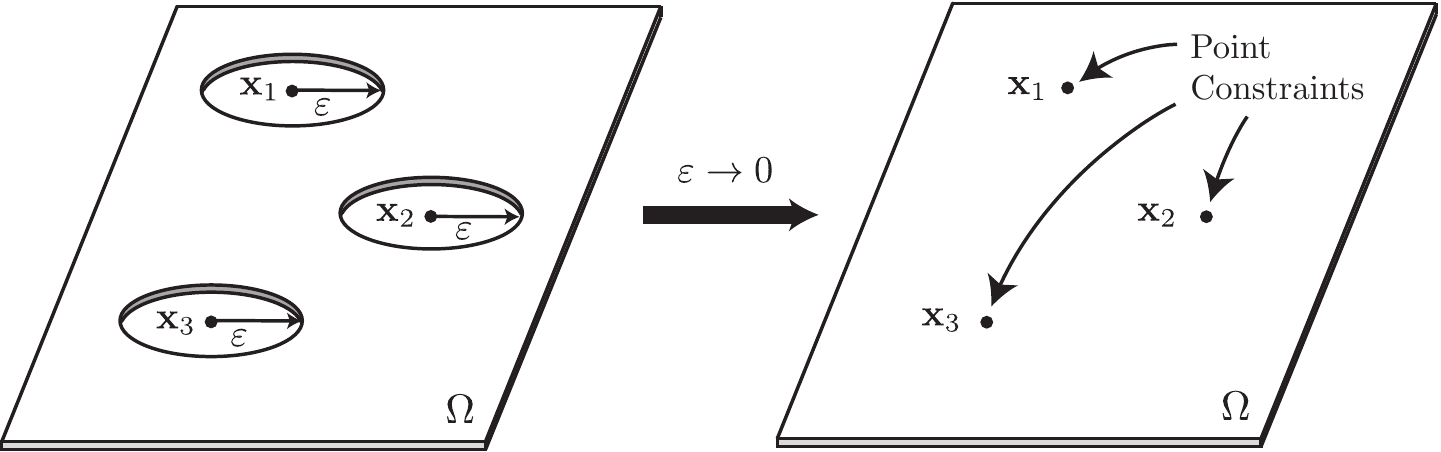}
\parbox{0.75\textwidth}{\caption{In the limit of vanishing hole radius $\eps\to0$, a point constraint $u(\xx_j)=0$ must be enforced at each of the hole centers for $j = 1,\ldots, M$.\label{fig:IntroShrink}}}
\end{figure}

In the present work, we consider a finite collection of $M$ defects or punctures on~\eqref{eqn:introA} with the conditions
\begin{equation}\label{eqn:introC} 
u(\xx_j) = 0, \qquad j = 1,\ldots, M.
\end{equation}
\esub
These \emph{point constraints} arise in singular perturbation studies
of~\eqref{eqn:introA} in the presence of $M$ small circular
perforations of radius $\eps$ (cf.~Fig.~\ref{fig:IntroShrink}). As the
radius $\eps$ of the perforations shrink to zero, the behavior of the
limiting eigenvalue $\lambda_{\eps}$ as
$\eps\to0$ satisfies~\cite{KLW,LWK,LHS,CN01}
\begin{equation}\label{HoleBehavior}
  \lambda_{\eps} = \lambda + 4\pi \nu \sum_{j=1}^M 
  |\nabla u(\xx_j)|^2 + \bigoh(\nu^2), \qquad \nu = -\frac{1}{\log\eps},
\end{equation}
where $(\lambda,u)$ satisfies (\ref{eqn:introA}-\ref{eqn:introB}) plus the point
constraints~\eqref{eqn:introC}. In the degenerate case $\sum_{j=1}^M
|\nabla u(\xx_j)|^2 = 0$, equation~\eqref{HoleBehavior} is not valid
and a separate limiting form can be derived~\cite{CN01,KLW}.  The fact
that the clamping condition on each perforation leaves an imprint as
the radius shrinks to zero (Fig.~\ref{fig:IntroShrink}) implies that no
matter how small a perforation is, the vibrational characteristics are
distinct from the no hole problem
\begin{equation}\label{HoleFree}
  \Delta^2 u^{\star} = \lambda^{\star} u^{\star}, 
  \quad \xx \in\Omega; 
  \qquad u^{\ast} = \partial_\nn u^{\star} = 0, 
  \quad \xx \in\partial\Omega; 
  \qquad \int_{\Omega} {u^{\star}}^2 d\xx = 1.
\end{equation}

\begin{figure}[htbp]
\centering
\includegraphics[width = 0.45\textwidth]{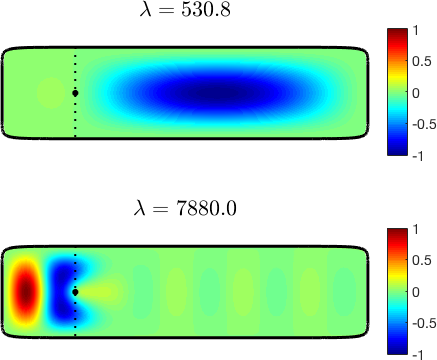}
\parbox{0.75\textwidth}{ \caption{The localization of two eigenfunctions
by a single clamped point, located at the black point, in a rectangular
domain. In each case, the eigenfunction is essentially zero on one side
of the clamping point. See~\cite{FM} and Sec.~\ref{sec:rectangle} for
more details.  \label{fig:introClamp} }}
\end{figure}

%
%
%
 %The specific name given to these localized structures is application dependent, however, in the rest of this article we refer to these localized features as punctures. 
%
The discontinuous limiting behavior of~\eqref{HoleBehavior} is
qualitatively different from the spectral problem for the Laplacian in
the presence of small perturbing holes~\cite{F,KTW,O,WHK,WK}. A
consequence of the point constraints~\eqref{eqn:introC} is that the
eigenfunctions $u(\xx)$ are not necessarily smooth but satisfy local
conditions
\begin{equation}\label{behaviorLocal}
  u(\xx) = \alpha_j |\xx - \xx_j|^2 \log|\xx - \xx_j| + 
    \bigoh(1), \qquad \xx\to\xx_j; \qquad j = 1,\ldots,M,
\end{equation}
where the constants $\{\alpha_j\}_{j=1}^M$ reflect the strength of each
puncture and depend on the domain $\Omega$ and the clamping locations
$\{\xx_j\}_{j=1}^M$. The difference between the punctured eigenvalues
$\lambda$ of~\eqref{eqn:intro} and the puncture free eigenvalues
$\lambda^{\star}$ of~\eqref{HoleFree} satisfies (cf.~\cite{LHS}) 
\begin{align}
  \label{EigDifference}
  (\lambda - \lambda^{\star}) \braket{u,u^{\star}} = 
  -8\pi\sum_{j=1}^M \alpha_j u^{\star}(\xx_j), \qquad 
  \braket{u,u^{\star}} = \int_{\Omega} u(\xx) u^{\star}(\xx)\, d\xx.
\end{align}
The presence of clamped locations also has a profound localizing effect
on the eigenfunctions. In a rectangular domain with a single clamped
point located along the long axis, the effect of clamping
on~\eqref{eqn:intro} has been observed (cf.~\cite{FM}) to partition
$\Omega$ into two distinct domains on the left and right of the clamping
location, as shown in Fig~\ref{fig:introClamp}. One aim of this work is
to numerically investigate the global effects that point constraints
have on the eigenfunctions of~\eqref{eqn:intro} in a variety of
different planar geometries.
 
Fourth-order eigenvalue problems~(Equations~\eqref{eqn:intro} and
\eqref{HoleFree}) exhibit other qualitatively different properties
compared to the well-understood Laplacian counterpart. For example, the
fundamental eigenfunction of~\eqref{eqn:intro}, ie.~the mode associated
with the lowest eigenvalue, is not necessarily single signed~\cite{CD1,
CD2, S, Coffman82, Gazzola2010, Grunau2014}. In contrast, the
fundamental eigenfunction of the Laplacian is always single signed and
the corresponding eigenvalue is simple~\cite{Evans2010,Gilbarg1998}.
An elementary example of this phenomenon is the annular domain
$\eps<r<1$ in which the radially symmetric and mode $1$ eigenvalues of
the bi-Laplacian cross at $\eps^{-1}\approx762.36$~\cite{CD2}.
Correspondingly, for $\eps^{-1}>762.36$, the fundamental eigenfunction
has multiplicity two and one nodal line. Also, in domains with a
corner, the first eigenfunction may possess an infinite number of nodal
lines~\cite{Coffman82}. Many numerical methods have been developed to
treat fourth-order eigenvalue problems in view of these
characteristics~\cite{brown2000, AA, CD, LC2010,
jia-kro-qua2013,Zhao2002}.

The main goal of this paper is to introduce a novel high-order boundary
integral equation method for the numerical solution of~\eqref{eqn:intro}
in the presence of a finite collection of punctures \eqref{eqn:introC}.
High-order methods for computing eigenvalues of the Laplacian and
Helmholtz equations in two and three dimensions have been developed with
domain decomposition methods~\cite{bet2007, des-tol1983, dri1997},
radial basis functions~\cite{pla-dri2004}, boundary integral
equations~\cite{bac2003, ste-ung2009, dur-ned-oss2009}, the method of
particular solution~\cite{bar2009, fox-hen-mol1967, kar2001}, the
Dirichlet to Neumann map~\cite{bar-has2014}, and chebfun~\cite{drum}.
The method of fundamental solutions has also been used to compute
eigenvalues of the biharmonic equation~\cite{mar-les2005, AA}.  However,
none of these works consider the eigenvalue problem with clamped points.
We extend the work of one of the previous authors~\cite{LHS} where a
finite difference method coupled with an inexact Newton method is used
to solve~\eqref{eqn:intro} in the unit circle with symmetrically chosen
clamped points.  Owing to the accuracy and robustness of the boundary
integral equation methods, our new method forms third-order solutions
of~\eqref{eqn:intro} in smooth two-dimensional geometries, including
multiply-connected geometries (Figure~\ref{fig:couetteModes}), and with
a large assortment of clamping locations.  

Using our new method, we demonstrate the dependence of $\lambda$ on the
number $M$ and locations $\{\xx_1, \ldots,\xx_M \}$ of the puncture
sites for a variety of planar regions $\Omega \subset \mathbb{R}^2$.
In particular, we investigate two effects that clamped points have on
the vibrational properties of plates with various regular and irregular
geometries. Our first observation is that by specific location of
punctures, the vibrational properties can be dramatically altered---in
particular, undesirable frequencies of vibration can be tuned out by
deliberate location of clamped points at nodal lines of the unclamped
eigenfunction $u^{\star}$ of \eqref{HoleFree}. Our second observation,
extending previous results in~\cite{FM} for rectangular domains, is
that mode confinement occurs in a variety of two dimensional geometries.

The outline of the paper is as follows. In Section~\ref{sec:methods} we
describe the details of a boundary integral method for
solving~\eqref{eqn:intro}. In Section~\ref{sec:algorithms}, the
implementation details are discussed and third-order convergence of the
method is verified for a closed-form solution of~\eqref{eqn:intro}.  In
Section~\ref{sec:numerics}, we apply our method to a disk, rectangles,
an ellipse, a non-symmetric shape, and a multiply-connected region.
Finally, in Section~\ref{sec:conclusions} we discuss the results and
areas of future investigations.

%%%%%%%%%%%%%%%%%%%%%%%%%%%%%%%%%%%%%%%%%%%%%%%%%%%%%%%%%%%%%%%%%%%%%%%%
\section{Integral equation formulation of the clamped eigenvalue
problem}
\label{sec:methods}
%%%%%%%%%%%%%%%%%%%%%%%%%%%%%%%%%%%%%%%%%%%%%%%%%%%%%%%%%%%%%%%%%%%%%%%%
In this section, we first compute and analyze the fundamental solution
of the modified biharmonic operator $\Delta^2 - \lambda$.  We then use
the fundamental solution to reformulate equation~\eqref{eqn:intro} as a
system of second-kind boundary integral equations with compact integral
operators.

%%%%%%%%%%%%%%%%%%%%%%%%%%%%%%%%%%%%%%%%%%%%%%%%%%%%%%%%%%%%%%%%%%%%%%%%
\subsection{Fundamental solution}
\label{sec:fundSoln}
We require the fundamental solution $G(\xx,\yy)$ of the modified biharmonic operator satisfying
\begin{align*}
  \Delta^2 G - \mu^4 G = \delta(\xx-\yy), \qquad \xx \in \mathbb{R}^2,
\end{align*}
where $\lambda = \mu^4$.  The factorization $\Delta^2 - \mu^4 =
(\Delta - \mu^2)(\Delta + \mu^2)$, and the fact the fundamental
solution is radially symmetric, imposes that $G(\xx,\yy)$ is a linear
combination of the Bessel functions $J_0(\mu \rho)$, $Y_0(\mu \rho)$,
$I_0(\mu \rho)$, and $K_0(\mu \rho)$, where $\rho = |\xx-\yy|$.  Using
a linear combination of the two singular Bessel functions that decay as
$r \rightarrow \infty$, the fundamental solution centered at $\yy$ is
of the form
\begin{align*}
  G(\xx,\yy) = c_{1}Y_{0}(\mu|\xx-\yy|) + c_{2}K_{0}(\mu|\xx-\yy|).
\end{align*} 
To find the appropriate constants $c_1, c_2$, we use the identities
$(\Delta + \mu^{2})Y_{0}(\mu|\xx-\yy|) = -4\delta(\xx-\yy)$ and $(\Delta
- \mu^{2})K_{0}(\mu|\xx - \yy|)  = -2\pi\delta(\xx-\yy)$, and compute
the fundamental solution by solving
\begin{align*}
  (\Delta - \mu^{2})(\Delta + \mu^{2})c_{1}Y_{0}(\mu|\xx-\yy|) + 
  (\Delta + \mu^{2})(\Delta - \mu^{2})c_{2}K_{0}(\mu|\xx-\yy|)=
  \delta(\xx-\yy).
\end{align*}
This calculation reveals that the fundamental solution of $\Delta^{2} -
\mu^{4}$ centered at $\yy$ is 
\begin{align}
  \label{FundamentalGreens}
  G(\xx,\yy) = -\frac{1}{8\mu^{2}} Y_{0}(\mu|\xx-\yy|) - 
            \frac{1}{4\pi\mu^{2}} K_{0}(\mu|\xx-\yy|).
\end{align}
We will be using $G$ in an indirect integral equation formulation, and
this will require the behavior of the fundamental solution when $\xx
\rightarrow \yy$.  Without loss of generality, we take $\yy=\textbf{0}$
and expand the fundamental solution for small $|\xx|$.  Using small
argument approximations of the Bessel functions
(cf.~\cite{abr-ste1964}), we have
\begin{align*}
  G(\xx,\textbf{0}) =  
      \frac{|\xx|^2}{8\pi}\log |\xx| \left(1 +
      \mathcal{O}(|\xx|^{4})\right)
      +\frac{|\xx|^{2}}{8\pi}\left(-1 + \gamma +
      \log\left(\frac{\mu}{2}\right) + \mathcal{O}(|\xx|^{4})\right),
      \quad \mbox{as} \quad |\xx|\to0,
\end{align*}
where $\gamma \approx 0.5772156649$ is Euler's constant. As mentioned
in the introduction, a key behavior of the solution
of~\eqref{eqn:intro} is the local behavior~\eqref{behaviorLocal} near
each of the defects.  Since the fundamental solution satisfies this
required behavior, the solution of~\eqref{eqn:intro} can be written as
\begin{align}\label{eqn:LinSep}
  u(\xx) = u_S(\xx)+ u_R(\xx), \qquad u_S(\xx) = 8\pi\sum_{j=1}^M \alpha_j G(\xx,\xx_j),
\end{align}
where $G(\xx,\yy)$ is given in~\eqref{FundamentalGreens}. In Section
\ref{sec:Newton}, we describe an inexact Newton method to find the
strength of the defects $\{\alpha_{j}\}_{j=1}^M$ and the eigenvalues
$\lambda$.  The decomposition~\eqref{eqn:LinSep} of the solution as the
sum of a singular and regular part allows for the local
behavior~\eqref{behaviorLocal} to be precisely enforced while the
regular part $u_R$ satisfies the homogeneous fourth-order PDE
\bsub\label{eqn:regularPDE}
\begin{gather}
  \label{eqn:regularPDE_A}
  \Delta^{2}u_{R} - \lambda u_{R} = 0, \quad \xx \in \Omega; \\[5pt]
  \label{eqn:regularPDE_B}  u_{R} = -u_{S}, \qquad \partial_{n} u_{R} = -\partial_{n} u_{S}, 
  \quad \xx \in \partial \Omega,
\end{gather}
\esub
where $u_S$ is specified in \eqref{eqn:LinSep}. We note that in~\cite{LHS}, the singular part was chosen to be
\begin{align*}
  u_{S}(\xx) = \sum_{j=1}^{M} \alpha_{j} |\xx - \xx_{j}
    |^2\log |\xx - \xx_{j}|.
\end{align*}
While this choice has the correct local behavior~\eqref{behaviorLocal},
it leads to a forcing term in the PDE for $u_R$ that, for a boundary
integral equation method, is prohibitive. However, the boundary
conditions~\eqref{eqn:regularPDE_B} in our new formulation depends
nonlinearly on the unknown eigenvalue $\lambda$.

Once the functions $u_S$ and $u_R$ are computed, they can be easily
evaluated at the locations of the clamped points.  This is used to
iteratively solve the non-linear equation (Section~\ref{sec:Newton})
\begin{align}
\label{eqn:mismatch}
  F(\mathbf{z}) = \left[
  \begin{array}{c}
    u_{S}(\xx_{1}) + u_{R}(\xx_{1}) \\
    \vdots \\
    u_{S}(\xx_{M}) + u_{R}(\xx_{M}) \\
    \alpha_{1}^{2} + \cdots + \alpha_{M}^{2} - 1
  \end{array}
  \right] = \left[
  \begin{array}{c}
    0 \\ \vdots \\ 0 \\ 0
  \end{array}
  \right],
\end{align}
where $\mathbf{z} = (\alpha_1,\ldots,\alpha_M,\lambda)$.  The particular
normalization condition $\sum_{j=1}^M \alpha_j^2 = 1$ is chosen purely
for ease of implementation. Once a solution is obtained, the
eigenfunction can be normalized according to~\eqref{eqn:introA} or any
other condition.

%%%%%%%%%%%%%%%%%%%%%%%%%%%%%%%%%%%%%%%%%%%%%%%%%%%%%%%%%%%%%%%%%%%%%%%%
\subsection{Computing the regular solution $u_R$}
\label{sec:layer_pots}
%%%%%%%%%%%%%%%%%%%%%%%%%%%%%%%%%%%%%%%%%%%%%%%%%%%%%%%%%%%%%%%%%%%%%%%%
Equation~\eqref{eqn:regularPDE} is linear and homogeneous, so it can be
recast in terms of a boundary integral equation.  In this section, we
describe appropriate layer potentials. Since the PDE is fourth-order, a
sum of two linearly independent layer potentials must be used.  The
regular part $u_R$ is written as
\begin{align}
  u_{R}(\xx) = \int_{\bd\Omega} G_{1}(\xx,\yy)\sigma_{1}(\yy) ds_{\yy} +
               \int_{\bd\Omega} G_{2}(\xx,\yy)\sigma_{2}(\yy) ds_{\yy},
  \label{eqn:layerPot}
\end{align}
where $G_{1}$ and $G_{2}$ are linear combinations of $G$ and its partial
derivatives.  The choice of $G_{1}$ and $G_{2}$ determines the nature of
the boundary integral equation which plays a crucial role on the
conditioning of the linear system that arises after discretization.  In
particular, $G_{1}$ and $G_{2}$ should be chosen so that the resulting
boundary integral equation is of the second-kind with compact integral
operators.  This means that the limiting values of the layer potential
ansatz~\eqref{eqn:layerPot} must have jumps that are proportional to
$\sigma_1$ and $\sigma_2$ as $\xx \rightarrow \bd\Omega$, and the
kernels must be integrable.

To find kernels $G_1$ and $G_2$ with these desired results, we use the
work of Farkas~\cite{far1989} who formulated the desired second-kind
integral equations for the fourth-order biharmonic equation.  For the
biharmonic equation with Dirichlet and Neumann boundary conditions,
Farkas proposed the kernels
\begin{align*}
  G_{1}(\xx,\yy) &= G_{\nn\nn\nn} + 3G_{\nn\ttau\ttau}, \\
  G_{2}(\xx,\yy) &= \Delta G - 2G_{\nn\nn},
\end{align*}
where the normal vector $\nn$ and tangent vector $\ttau$ are taken with
respect to the source point $\yy$.  Since the leading order singularity
of $G$, $\frac{1}{8\pi}|\xx|^{2}\log|\xx|$, is equal to the fundamental
solution of the two-dimensional biharmonic equation, the jumps in the
layer potential~\eqref{eqn:layerPot} agree, to a first approximation,
with the jumps found by Farkas.  In particular, any additional jumps in
$G_1$ and $G_2$ will result from the higher-order terms in the expansion
of $G$.  Since the higher-order terms contain singularities of strength
no less than $|\xx|^{6}\log|\xx|$, no additional jumps will be present
as long as $G_1$ and $G_2$ do not involve derivatives of order six or
higher.  Since the derivatives $G_1$ and $G_2$ are no more than
third-order, the jumps of $G_1$ and $G_2$ will agree with those computed
by Farkas.

%%%%%%%%%%%%%%%%%%%%%%%%%%%%%%%%%%%%%%%%%%%%%%%%%%%%%%%%%%%%%%%%%%%%%%%%
\subsection{Explicit expressions of the kernels}
\label{sec:kernels}
%%%%%%%%%%%%%%%%%%%%%%%%%%%%%%%%%%%%%%%%%%%%%%%%%%%%%%%%%%%%%%%%%%%%%%%%
For $\xx,\yy \in \bd\Omega$, we require the four kernels
\begin{align*}
  G_{11}(\xx,\yy) &= G_{1}(\xx,\yy), \\
  G_{12}(\xx,\yy) &= G_{2}(\xx,\yy), \\
  G_{21}(\xx,\yy) &= \pderiv{}{\nn_{\xx}}G_{1}(\xx,\yy), \\
  G_{22}(\xx,\yy) &= \pderiv{}{\nn_{\xx}}G_{2}(\xx,\yy).
\end{align*}
Substituting the fundamental solution~\eqref{FundamentalGreens} into
these expressions, and using the identities
\begin{align*}
  \pderiv{}{\nn}(\rdotn) = -1, \qquad
  \pderiv{}{\nn}\rho = -2 \frac{\rdotn}{\rho^{2}}, \qquad
  \pderiv{}{\nn_\xx}(\rdotn) = 1, \qquad
  \pderiv{}{\nn_\xx}\rho = -2 \frac{\rdotnx}{\rho^{2}},
\end{align*}
where $\rr = \xx - \yy$, $\rho = |\rr|$, and similar identities for
the tangential derivatives, the kernels $G_{11}$ and $G_{12}$ are
\begin{align*}
  G_{11} = -\frac{1}{4\pi\mu^2} &\left(
    3\mu^{3}K_{1}(\mu\rho)\frac{(\rdotn)}{\rho} - 
    2\mu^{3}K_{1}(\mu\rho)\frac{(\rdotn)^{3}}{\rho^3} +
    6\mu^{2} K_{0}(\mu\rho)\frac{(\rdotn)}{\rho^{2}} \right. \\
    &\left.
    -8\mu^{2} K_{0}(\mu\rho)\frac{(\rdotn)^{3}}{\rho^{4}} - 
    16\mu K_{1}(\mu\rho)\frac{(\rdotn)^{3}}{\rho^{5}} + 
    12\mu K_{1}(\mu\rho)\frac{(\rdotn)}{\rho^{3}}
    \right) \\
    -\frac{1}{8\mu^{2}}&\left(
    -3\mu^{3}Y_{1}(\mu\rho)\frac{(\rdotn)}{\rho} + 
    2\mu^{3}Y_{1}(\mu\rho)\frac{(\rdotn)^{3}}{\rho^3} -
    6\mu^{2} Y_{0}(\mu\rho)\frac{(\rdotn)}{\rho^{2}} \right. \\
    &\left.
    +8\mu^{2} Y_{0}(\mu\rho)\frac{(\rdotn)^{3}}{\rho^{4}} - 
    16\mu Y_{1}(\mu\rho)\frac{(\rdotn)^{3}}{\rho^{5}} + 
    12\mu Y_{1}(\mu\rho)\frac{(\rdotn)}{\rho^{3}}
    \right), \\
  G_{12} = -\frac{1}{4\pi}&\left(1 - 2\frac{(\rdotn)^2}{\rho^{2}}\right)
    \left(K_{0}(\mu\rho) + \frac{2}{\mu\rho}K_{1}(\mu\rho)\right) 
    +\frac{1}{8}\left(1 - 2\frac{(\rdotn)^{2}}{\rho^{2}}\right)
    \left(Y_{0}(\mu\rho) - \frac{2}{\mu\rho}Y_{1}(\mu\rho)\right).
\end{align*}
The expressions for $G_{21}$ and $G_{22}$ require one additional
derivative of $G_{11}$ and $G_{12}$.  For completeness, these lengthy expressions are given in Appendix~\ref{app:kernels}.

%%%%%%%%%%%%%%%%%%%%%%%%%%%%%%%%%%%%%%%%%%%%%%%%%%%%%%%%%%%%%%%%%%%%%%%%
\subsection{The boundary integral equation}
\label{sec:BIE}
%%%%%%%%%%%%%%%%%%%%%%%%%%%%%%%%%%%%%%%%%%%%%%%%%%%%%%%%%%%%%%%%%%%%%%%%

As discussed in Section~\ref{sec:fundSoln}, all four kernels $G_{ij}$
have the same asymptotic behavior as the fundamental solution of the
biharmonic equation.  Therefore, the boundary integral equation for
$\ssigma$ is identical to the boundary integral equation for the
biharmonic equation~\cite{far1989},
\begin{align}
  \label{eqn:bie}
  D(\xx)\ssigma(\xx) + \int_{\bd\Omega} A(\xx,\yy)\ssigma(\yy) ds_{\yy}
    =\gg(\xx),
\end{align}
where
\begin{align*}
  D(\xx) = \left(
  \begin{array}{cc}
    \displaystyle\frac{1}{2} & 0 \\ -\kappa(\xx) & \displaystyle\frac{1}{2}
  \end{array}
  \right),
\end{align*}
$\kappa(\xx)$ is the curvature of $\bd\Omega$ at $\xx$, and 
\begin{align*}
  \gg = -\left(
  \begin{array}{c}
    u_S \\ \partial_{\nn}u_S
  \end{array}
  \right),
  \quad 
  \ssigma = \left(
  \begin{array}{c}
    \sigma_{1} \\ \sigma_{2}
  \end{array}
  \right),
  \quad
  A = \left(
  \begin{array}{cc}
    G_{11} & G_{12} \\ 
    G_{21} & G_{22}
  \end{array}
  \right).
\end{align*}
To apply quadrature formulae, the limiting values of $G_{ij}$ as $\xx\to\yy$ are
required.  These can be found by applying L'H\^{o}pital's rule to each
of the four kernels.  For $\xx$, $\yy$ on $\bd\Omega$ we have
\begin{equation}
  \label{eqn:limits}
  \begin{aligned}
    \lim_{\yy \rightarrow \xx} G_{11}(\xx,\yy) &= 0, \\
    \lim_{\yy \rightarrow \xx} G_{12}(\xx,\yy) &=
    \frac{1}{4\pi}\kappa(\xx), \\
    \lim_{\yy \rightarrow \xx} G_{21}(\xx,\yy) &=
    -\frac{3}{4\pi}\kappa(\xx)^2, \\
    \lim_{\yy \rightarrow \xx} G_{22}(\xx,\yy) &=
    \frac{1}{2\pi}\kappa(\xx).
  \end{aligned}
\end{equation}

%%%%%%%%%%%%%%%%%%%%%%%%%%%%%%%%%%%%%%%%%%%%%%%%%%%%%%%%%%%%%%%%%%%%%%%%
\section{Numerical Methods}
\label{sec:algorithms}
%%%%%%%%%%%%%%%%%%%%%%%%%%%%%%%%%%%%%%%%%%%%%%%%%%%%%%%%%%%%%%%%%%%%%%%%
Here we describe a numerical method for solving the boundary integral
equation~\eqref{eqn:bie} (Section~\ref{sec:bie_dis}), applying an
inexact Newton method for~\eqref{eqn:mismatch}
(Section~\ref{sec:Newton}), and an algorithm for tracing the first
eigenvalue, $\lambda$, as clamped points are smoothly moved through the
geometry $\Omega$ (Section~\ref{sec:initial_guess}). 

%%%%%%%%%%%%%%%%%%%%%%%%%%%%%%%%%%%%%%%%%%%%%%%%%%%%%%%%%%%%%%%%%%%%%%%%
\subsection{Discretization of the integral equation}
\label{sec:bie_dis}
%%%%%%%%%%%%%%%%%%%%%%%%%%%%%%%%%%%%%%%%%%%%%%%%%%%%%%%%%%%%%%%%%%%%%%%%
We apply a standard collocation method to solve the second-kind boundary
integral equation~\eqref{eqn:bie}.  The boundary, $\bd\Omega$, is first
discretized at collocation points $\xx_{i}$,
$i=1,\ldots,N$.  To satisfy the boundary integral equation at these
collocation points, we require
\begin{align}\label{eq:integralMain}
  D(\xx_i)\ssigma(\xx_i) + \int_{\bd\Omega} A(\xx_i,\yy)\ssigma(\yy) ds_{\yy}
    =\gg(\xx_i).
\end{align}
The integral in~\eqref{eq:integralMain} is approximated with the
trapezoid rule where the abscissae are the collocation points which
yields the dense linear system
\begin{align*}
  D(\xx_{i}) \ssigma_i + \sum_{j=1}^{N} A(\xx_{i},\xx_{j})\Delta s_{j}
    \ssigma_{j} = \gg_{i},
\end{align*}
where $\ssigma_i = \ssigma(\xx_{i})$, $\gg_{i} = \gg(\xx_{i})$, and
$\Delta s_{j}$ is the Jacobian of the curve at point $\xx_{j}$.  The
limiting values from~\eqref{eqn:limits} are used for the diagonal terms
$A(\xx_i,\xx_i)$ of the linear system.

The convergence order of the method depends on the regularity of the
kernels $G_{ij}$.  The regularity of the kernels can be computed by
taking a simple geometry, such as the unit circle, fixing $\xx$, and
computing the limit as $\yy \rightarrow \xx$ of $G_{ij}(\xx,\yy)$ and
its derivatives.  These calculations reveal that
\begin{align*}
  G_{11} \in C^3, \qquad G_{12} \in C^3, \qquad G_{21} \in C^1, \qquad
  G_{22} \in C^3.
\end{align*}
The accuracy of the trapezoid rule for a periodic $C^k$ function is
$k+2$, so we expect third-order convergence because of the $C^1$
regularity of $G_{21}$.  Higher-order accuracy can be achieved by using
specialized quadrature~\cite{alp1999, kap-rok1997} designed for
functions with weak logarithmic singularities. Once values for the
density function $\ssigma_j$ are computed, we can compute $u_R(\xx)$ for
any $\xx \in \Omega$ with spectral accuracy.  In particular, we compute
the value at the clamped locations with the trapezoid rule to yield that
\begin{align*}
  u_{R}(\xx) &= \int_{\bd\Omega} G_{1}(\xx,\yy)\sigma_{1}(\yy) ds_{\yy} +
             \int_{\bd\Omega} G_{2}(\xx,\yy)\sigma_{2}(\yy) ds_{\yy}\\[5pt]
             &\approx \frac{2\pi}{N} \sum_{j} \left( G_{1}(\xx,\yy_j) \sigma_{1_j} +
                G_{2}(\xx,\yy_j) \sigma_{2_j}
                \right)\Delta s_j.
\end{align*}
If a target point $\xx$ is sufficiently close to $\bd\Omega$, then the
accuracy of the trapezoid rule will be diminished due to large derivatives
in $G_{i}(\xx,\yy)$.  In this case, we simply upsample the geometry and
density functions so that sufficient accuracy can be achieved at the
clamped locations.

%%%%%%%%%%%%%%%%%%%%%%%%%%%%%%%%%%%%%%%%%%%%%%%%%%%%%%%%%%%%%%%%%%%%%%%%
\subsection{Nonlinear solvers}\label{sec:Newton}
To solve the nonlinear equation~\eqref{eqn:mismatch} for
$\{\alpha_j\}_{j=1}^{M}$ and $\lambda$, we apply one of two strategies.
First, in symmetric cases such as the disk geometry, if the clamped
points are equidistributed in the azimuthal direction at a fixed radius,
then $\alpha_1 = \cdots = \alpha_M$.  Therefore, $\alpha_j =
M^{-\frac{1}{2}}$ for $j = 1,\ldots,M$, and the only parameter remaining
is $\lambda$. For such a case, and any scenario in which symmetry
considerations reduce the unknown to just $\lambda$, a bisection method
can be applied to reliably solve~\eqref{eqn:mismatch} since convergence
to the desired root is guaranteed for an appropriately chosen initial
interval. This method is generally preferred in cases where all the
$\alpha_j$ are equal and the single unknown is the eigenvalue itself.

Second, when symmetry can not be assumed, we apply an inexact Newton's
method to~\eqref{eqn:mismatch}.  In our calculations, the Jacobian
matrix $J$ of $F$ is formed by finite difference approximations which we
have found to be accurate and efficient.
%\begin{equation}\label{eqn:NonlinearSys}
%  F(\textbf{z}) = 
%  \begin{bmatrix}
%    u_{S}(\xx_{1}) + u_{R}(\xx_{1}) \\
%    \vdots \\
%    u_{S}(\xx_{M}) + u_{R}(\xx_{M}) \\
%    \alpha_{1}^{2} + \cdots + \alpha_{M}^{2} - 1
%  \end{bmatrix}
%  = \textbf{0},
%\end{equation}
%for a solution vector $\textbf{z} = (\alpha_1,\ldots,\alpha_M,\lambda)$.

\begin{figure}[htbp]
\centering
\begin{tikzpicture}[scale=1.0]

\begin{axis}[
  title = {},
  xmin = 14,
  xmax = 2300,
  xtick = {32,128,512,2048},
  xticklabels = {32,128,512,2048},
  xmode = log,
  xlabel = {\Large $N$},
  ymin = 1.0e-9,
  ymax = 1.0e-2,
  ytick = {1e-8,1e-6,1e-4,1e-2},
  ymode = log,
  ylabel = { \Large $\mathcal{E}_{\textrm{rel}}$},
  ylabel style = {yshift = 1pt},
  ylabel style = {xshift = 11pt},
  ylabel style = {rotate=-90},
  xlabel style = {xshift = 12pt},
  ]

% relative error
\addplot [color=black,mark=*,solid,line width=2] coordinates{
(16,6.5562e-03)
(32,8.2126e-04)
(64,1.0302e-04)
(128,1.2893e-05)
(256,1.6121e-06)
(512,2.0153e-07)
(1024,2.5191e-08)
(2048,3.0763e-09)
};

% line of slope 3
\addplot [color=red,dashed,line width=2] coordinates{
(32,3.0000e-3)
(1024,9.1553e-8)
};

\end{axis}

\end{tikzpicture}
\parbox{0.75\textwidth}{\caption{The relative error (black) of our
numerical method when using the bisection method to find the first
eigenvalue of \eqref{eqn:intro} with a single clamped point at the
center of the unit disk.  A line of slope $-3$ (red) indicates the
expected third-order convergence.\label{fig:ErrorDisk}}}
\end{figure}
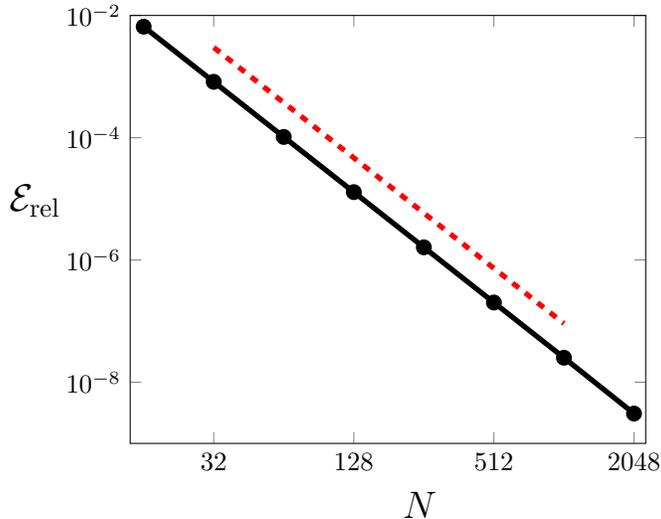

We validate the method with the unit disk geometry. A closed-form
solution of~\eqref{eqn:intro} can be developed in the special case $M=1$
and $\xx_1 =(0,0)$. In a similar manner to the construction of the
fundamental solution~\eqref{FundamentalGreens}, a linear combination of
$K_0$ and $Y_0$ can be chosen to eliminate the logarithmic singularity
at the origin.  Therefore radially symmetric eigenfunctions
of~\eqref{eqn:intro} are a combination of $Y_0(\mu \rho),$ $K_0(\mu \rho)$, 
$J_0(\mu \rho)$ and $I_0(\mu \rho)$ with $\rho = |\xx|$. The eigenfunctions that are finite at the
origin and satisfy $u(0)=0$ and $u(1) = \partial_\rho u(1)=0$ are
\begin{equation*}
u(\rho) = A \left[J_0(\mu_{0,n} \rho) - I_0(\mu_{0,n} \rho) -
\left(\frac{J_0(\mu_{0,n}) - I_0(\mu_{0,n})}{\frac{2}{\pi} K_0(\mu_{0,n})
+ Y_0(\mu_{0,n})} \right) \left(\frac{2}{\pi} K_0(\mu_{0,n} \rho) +
Y_0(\mu_{0,n} \rho)\right) \right],
\end{equation*}
where $A$ is a normalization constant and the eigenvalues $\lambda_{0,n}
= \mu_{0,n}^4$ satisfy the relationship
\begin{equation}\label{disk_exact}
\big(J_0(\mu_{0,n}) - I_0(\mu_{0,n})\big) \left( \frac{2}{\pi}
K_1(\mu_{0,n}) + Y_1(\mu_{0,n}) \right) = \big(J_1(\mu_{0,n}) +
I_1(\mu_{0,n})\big) \left( \frac{2}{\pi} K_0(\mu_{0,n}) + Y_0(\mu_{0,n})
\right).
\end{equation}
The smallest positive solution of~\eqref{disk_exact} gives rise to the
eigenvalue $\lambda_{\textrm{true}} \approx 516.9609.$  This solution
provides a benchmark against which the efficacy of our numerical method
can be verified.  We compute the relative error
$\mathcal{E}_{\textrm{rel}}$ between the numerically determined value of
$\lambda_{\textrm{num}}$ and the exact value $\lambda_{\textrm{true}}$.
%\begin{equation*}
%  \mathcal{E}_{\textrm{rel}} = \left|\frac{\lambda_{\textrm{num}} - 
%  \lambda_{\textrm{true}}}{\lambda_{\textrm{true}}}\right|.
%\end{equation*}
In Fig.~\ref{fig:ErrorDisk}, the numerical error scales $\bigoh(N^{-3})$
as the number of boundary points $N$ increases which agrees with our
expected third-order convergence.  In this example, the bisection method
was used, and the strength of the singularity is $\alpha = 1$.

%%%%%%%%%%%%%%%%%%%%%%%%%%%%%%%%%%%%%%%%%%%%%%%%%%%%%%%%%%%%%%%%%%%%%%%%
\subsection{Initialization, parameterization of puncture patterns, and
arclength continuation}
\label{sec:initial_guess}

The solution of the nonlinear system~\eqref{eqn:mismatch} by Newton's
Method relies on good initial iterates. In addition, a careful selection
of the initial guess is necessary to reliably locate the lowest mode of
the punctured problem~\eqref{eqn:intro}. For the unit circle, we start
with the clamped points at the center of the circle and initialize
Newton iterations for~\eqref{eqn:mismatch} with the known eigenvalue
$\lambda \approx 516.9609$ for a single clamped point at the origin.
For other geometries, we start the clamped point near $\bd\Omega$.  In
this scenario, equation~\eqref{eqn:mismatch} is initialized with a mode
of the unclamped problem~\eqref{HoleFree} calculated from a low-accuracy
finite element approximation~\cite{KI78}.  Once a solution
of~\eqref{eqn:intro} has been generated, the punctures are gradually
moved, and~\eqref{eqn:mismatch} is repeatedly solved until the punctures
occupy a specified target set.

In the examples that follow, we compute eigenvalues $\lambda =
\lambda(r)$ of~\eqref{eqn:intro} for families of puncture patterns
described by a single parameter $r\geq0$. For reasons of efficiency
and to provide robustness to the Newton iterations, we use arc-length
adaptively to focus resolution at  sharp peaks of the curve as compared
to the surrounding areas. The algorithm to find points on the curve
$\lambda = \lambda(r)$ is initialized with a relatively large step size
$\mathrm{d}r$ with the concavity monitored until proximity to an extrema
is detected. Once an extrema of the curve is detected, $\mathrm{d}r$ is
reduced based on the current slope up to a minimum allowable step size.

%The step size is then a function of the current slope and the scale, along with the preset initial step size ($dy$) and
%minimum allowable step size ($rmin$) as shown in Equation \ref{dr}:
%\begin{equation} \label{dr}
%dr = \min(dy, \sqrt{r_{min}^2 + (dy*scale*slope^2});
%\end{equation}
%Because this is a numerical solution, there is some variation in the points around the curve. Thus, some memory was built into the code requiring multiple indications of concave down to enter the scaled dr mode and then multiple indications of concave up to leave it, greatly reducing the chances of unintended scaling or nonscaling of the step size. An additional fail safe was built into the code which compares a linear prediction of the lambda value to the calculated value at each step. It then reduces the time step and repeats the calculation if the difference of the values is greater than a user defined threshold.

%%%%%%%%%%%%%%%%%%%%%%%%%%%%%%%%%%%%%%%%%%%%%%%%%%%%%%%%%%%%%%%%%%%%%%%%
\section{Numerical Examples}
\label{sec:numerics}
%%%%%%%%%%%%%%%%%%%%%%%%%%%%%%%%%%%%%%%%%%%%%%%%%%%%%%%%%%%%%%%%%%%%%%%%

In this section we demonstrate the effectiveness of the method on
regular and irregular domains. To understand the role of clamping in the
eigenvalue problem, and interpret the results obtained with our
numerical method, we recall from~\eqref{EigDifference} that
\begin{align*}
%  \label{LinkLater}
 (\lambda - \lambda^{\star}) \braket{u,u^{\star}} = 
  -8\pi\sum_{j=1}^M \alpha_j u^{\star}(\xx_j), \qquad 
  \braket{u,u^{\star}} = \int_{\Omega} u(\xx) u^{\star}(\xx)\, d\xx,
\end{align*}
which relates the modes $(\lambda,u)$ of~\eqref{eqn:intro} to the
unclamped modes $(\lambda^{\star},u^{\star})$ of~\eqref{HoleFree}.  In
each of the examples that follow, we will use a low accuracy finite
element method~\cite{KI78} to obtain the required solutions
of~\eqref{HoleFree}. 

%This method utilizes an auxiliary variable
%$v^{\star} = \Delta u^{\star}$ so that~\eqref{HoleFree} can be written
%as a second-order system with boundary conditions $\partial_n u^{\star}
%= 0$ and $\partial_n v^{\star} = L u^{\star}$ on $\partial \Omega$.
%The regularizing parameter $L$ is chosen to be large and positive which
%enforces the condition $u^{\star}=0$ on $\partial\Omega$.

%\begin{equation}\label{eqn:EigsRegular}
%\begin{bmatrix} 
%\Delta & -1 \\ 0& \Delta
%\end{bmatrix}
%\begin{bmatrix} 
%u^{\star}\\ v^{\star}
%\end{bmatrix}
%=
% \lambda^{\star} \begin{bmatrix} 
%0 & 0\\ 1 & 0
%\end{bmatrix}
%\begin{bmatrix} 
%u^{\star}\\ v^{\star}
%\end{bmatrix}
%\quad \xx \in \Omega
%\qquad 
%\begin{matrix}
%\partial_n u^{\star} = 0\\
%\partial_n v^{\star} = \eta u^{\star}\\
%\end{matrix}
%\end{equation}

%%%%%%%%%%%%%%%%%%%%%%%%%%%%%%%%%%%%%%%%%%%%%%%%%%%%%%%%%%%%%%%%%%%%%%%%
\subsection{Unit circle}
The relationship~\eqref{EigDifference} shows how the distinct
eigenvalues and eigenfunctions of the clamped and unclamped
problems,~\eqref{eqn:intro} and~\eqref{HoleFree}, respectively, are
related. For each domain it is therefore important to consider the
solutions $(\lambda^{\star},u^{\star})$ to understand the effect of
puncture configurations.

\begin{figure}[htbp]
\centering
\includegraphics[width = 0.7\textwidth]{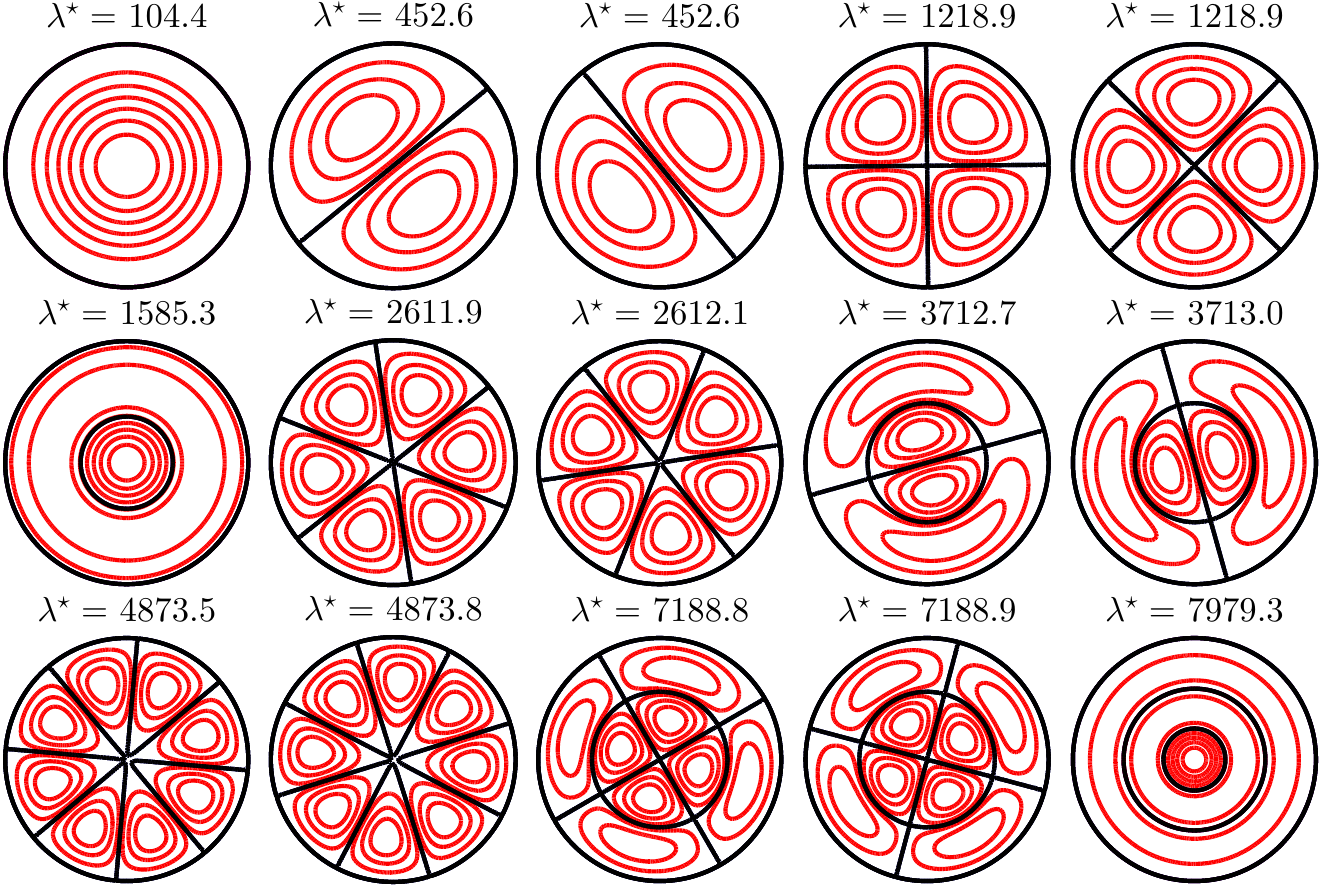}
\parbox{0.75\textwidth}{\caption{The contour lines of the first $15$
modes of the unclamped problem~\eqref{HoleFree} on the unit disk.  The
nodal lines ($u^{\star}=0$) are plotted in black.  Eigenfunctions are
repeated according to their multiplicity.  \label{fig:NoHoles}}}
\end{figure}

For the unit disk case, the solutions of problem \eqref{HoleFree} are
found by first factorizing $\Delta^2 - \mu^4 = (\Delta -\mu^2)(\Delta +
\mu^2)=0$ which indicates that the basis for the space of eigenfunctions
is
\begin{equation*}
%  \label{FullSpan}
 e^{im \theta} \{  J_m(\mu_{m,n} \rho), Y_m(\mu_{m,n} \rho), K_m(\mu_{m,n} \rho), I_m(\mu_{m,n} \rho) \}, \qquad \mu_{m,n} = \lambda_{m,n}^{1/4},
\end{equation*}
where $\rho = |\xx|$. The indices $m = 0, \pm 1, \pm2, \ldots$ indicate
the angular wavenumber (and number of angular nodal lines) where as
$n=0,1,2,\ldots$ counts the number of radial nodal lines for each
wavenumber. In the unclamped problem~\eqref{HoleFree}, the smooth
eigenfunctions satisfying $u^{\star} = \partial_{\rho} u^{\star} = 0$
on $\rho=1$ are
%\bsub\label{NoHolesExact}
\begin{align*}
%\label{NoHolesExact_a}
  u^{\star}_{m,n}(\rho,\theta) = e^{im\theta} \left[ J_m(\mu^{\star}_{m,n}\rho) - 
    \frac{J_m(\mu^{\star}_{m,n})}{I_m(\mu^{\star}_{m,n})} 
    I_m(\mu^{\star}_{m,n}\rho) \right],
\end{align*}
with the eigenvalues $\mu^{\star}_{m,n}$ determined by the relationship
\begin{align}\label{NoHolesExact_b}
J'_m(\mu^{\star}_{m,n})I_m(\mu^{\star}_{m,n}) = I'_m(\mu^{\star}_{m,n}) J_m(\mu^{\star}_{m,n}).
\end{align}
%\esub
The first four eigenvalues $\lambda^{\star}_{m,n} =
(\mu^{\star}_{m,n})^4$, found from the numerical solution
of~\eqref{NoHolesExact_b}, are
\begin{equation}\label{eigsDiskPFree}
\lambda^{\star}_{0,0} = 104.4, \qquad \lambda^{\star}_{1,0} = 452.0, \qquad \lambda^{\star}_{2,0} = 1216.4, \qquad \lambda^{\star}_{0,1} = 1581.7.
\end{equation}
In Fig.~\ref{fig:NoHoles}, the first few eigenfunctions are plotted with the nodal lines along which $u^{\star}=0$ highlighted.
For punctures away from the origin, we seek solutions
of~\eqref{eqn:intro} with parameterized puncture sets to minimize the
number of unknowns over which nonlinear iterations are processed. In the
disk case, our first example is a single ring of punctures given
explicitly as
\begin{align}
  \label{singleRing}
    \xx_j =  r \left( \cos\frac{2\pi j}{M}, 
    \sin \frac{2\pi j}{M} \right), \qquad j = 1,\ldots, M.
\end{align}

There is now a single parameter $r$ over which various configurations
can be investigated from $0<r<1$. From the eigenvalue $\lambda\approx
516.96$, which is the lowest positive solution of \eqref{disk_exact}, we
form an initial guess for the Newton iterations. Since our search
pattern~\eqref{singleRing} is radially symmetric, the puncture strengths
$\{\alpha_j\}_{j=1}^M$ can be assumed to be identical in this case. 
\begin{figure}[htbp]
\centering
{\subfigure[Single puncture ring pattern]{\includegraphics[width =
0.32\textwidth]{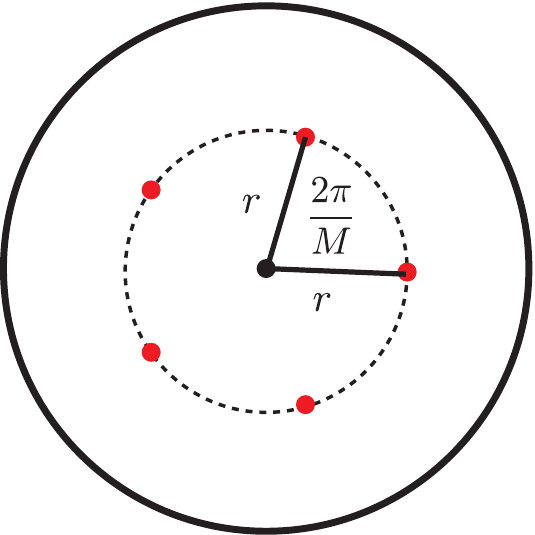}}\label{fig:SingleRing_A}}
\hspace{0.25in} {\subfigure[The lowest eigenvalue against ring
radius]{\includegraphics[width =
0.45\textwidth]{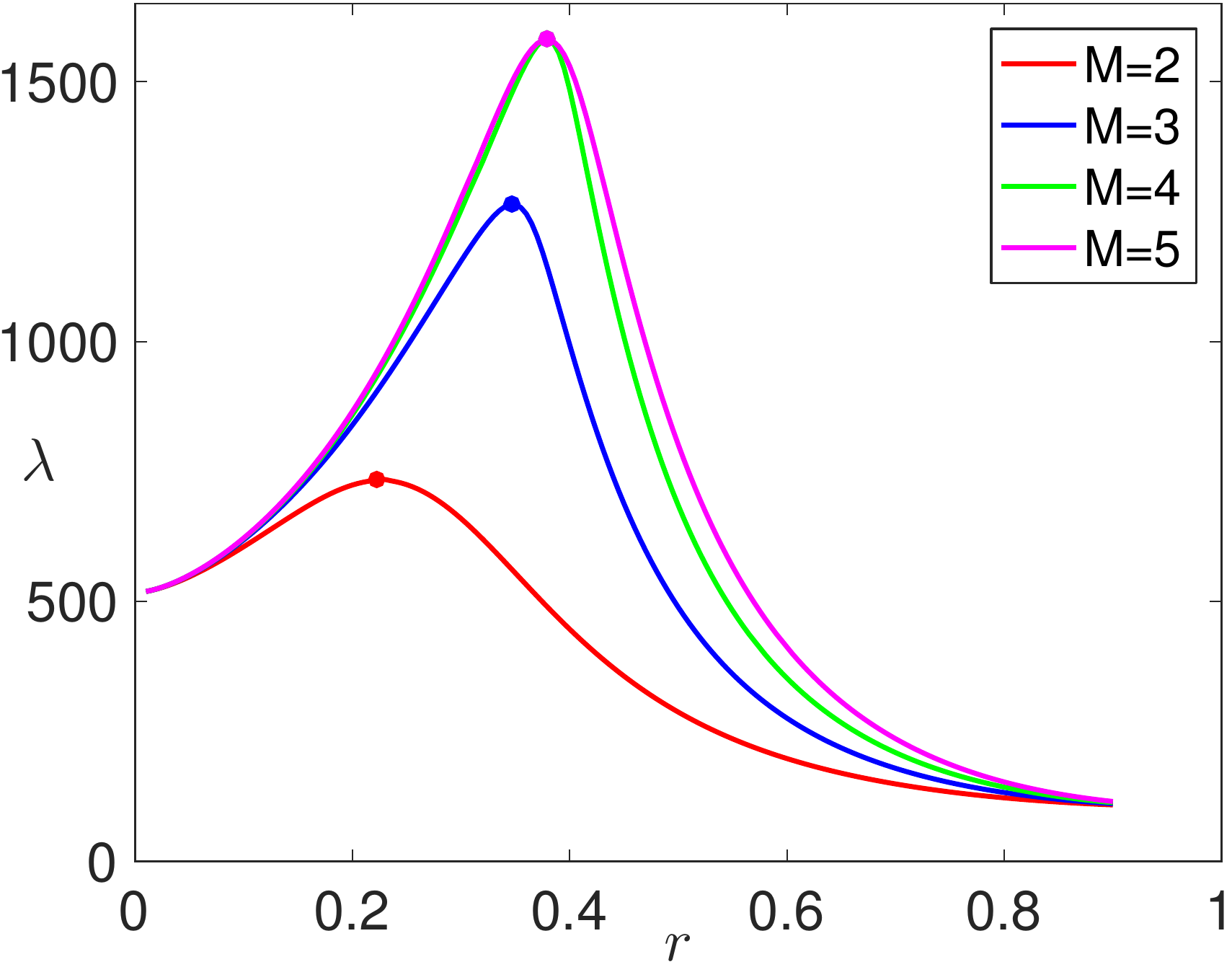}}\label{fig:SingleRing_B}}
\parbox{0.75\textwidth}{ \caption{Left: A disk with a single ring
puncture configuration (red dots). Right: The lowest radially symmetric
eigenvalue as a function of the puncture ring radius $r$.  The integral
equation~\eqref{eqn:bie} is discretized with $N=128$ boundary points.
\label{fig:SingleRing}} }
\end{figure}

In Fig.~\ref{fig:SingleRing} we see that as the eigenvalue is varied,
$\lambda$ attains a maximum value depending on the number of punctures
and the radius of the ring. In Table~\ref{tab:SingleRing}, we display
the maximum value of the lowest eigenvalue $\lambda_c$ and the critical
radius $r_c$ of puncture ring where it is attained. The results in
Table~\ref{tab:SingleRing} are in good agreement with~\cite{LHS} and
corroborate the observation that the maximum eigenvalue saturates as the
number of puncture points $M$ increases.

To explain the saturation effect, we recall
equation~\eqref{EigDifference} that relates the puncture and puncture
free modes. Consider a mode $(\lambda_k^{\star},u_k^{\star})$ of the
puncture free eigenfunction problem \eqref{HoleFree} such that
$u_k^{\star}(\xx_j) = 0$, for $j = 1,\ldots,M$. Then from
\eqref{EigDifference} we have that $(\lambda - \lambda_k^{\star})
\braket{u,u_k^{\star}} = 0$.
%
%\begin{equation*}
%  (\lambda - \lambda_k^{\star}) \braket{u,u_k^{\star}} = 0.
%\end{equation*}
This implies that either $\lambda = \lambda_k^{\star}$ or
$\braket{u,u_k^{\star}} = 0$ so that if $u\to u^{\star}_k$, then
$\lambda \to \lambda_k^{\star}$. From this we conclude that by centering
punctures on the nodal set of a puncture free eigenfunction
$u_k^{\star}$, then that eigenfunction becomes a mode of
\eqref{eqn:intro} thereby eliminating other modes form the spectrum.

\begin{table}[htbp]
\centering
\begin{tabular}{|c|ccccccc|}
\hline
$M$ & $2$ & $3$ & $4$ & $5$ & $6$ & $7$ & $8$\\
\hline
$r_c$ & $0.222$ & $0.348$ &  $0.379$&  $0.379$&  $0.379$&  $0.379$&  $0.379$\\
\hline
$\lambda_c$ & $734.96$ & $1264.2$ & $1581.5$ & $1581.5$ & $1581.5$ & $1581.5$ & $1581.5$\\
\hline
\end{tabular}
\parbox{0.75\textwidth}{ \caption{The maximum value $\lambda_c$ of the
lowest eigenvalue of the unit disk with a single ring of punctures.  The
values are obtained by numerical simulations with $N=128$ boundary
points.  \label{tab:SingleRing}} }
\end{table}

In this disk case with a single ring of punctures, saturation occurs
when the punctures are placed on the nodal set of
$u^{\star}_{0,1}(\mu^{\star}_{0,1} \rho)$---the mode with zero angular
nodal lines and one radial nodal line. From~\eqref{eigsDiskPFree}, we
see that the corresponding eigenvalue is $\lambda_{0,1}^{\star} =1581.7$
which agrees closely with the saturating value in Table~\ref{tab:SingleRing}.
The critical radius $r_c$ is found by solving
\begin{align*}
  u^{\star}_{0,1}(\mu^{\star}_{0,1} r_c) = 0 \quad 
    \Longrightarrow \quad r_c \approx 0.379,
\end{align*}
which agrees closely with the numerical results in
Table~\ref{tab:SingleRing}. 

The practical importance of this phenomenon is that undesirable
frequencies of the plate can be removed by specific placement of
clamping locations. For example, the placement of five clamped points
equally spaced on a circle of radius $r_c \approx 0.379$ yields that
the lowest mode of~\eqref{eqn:intro} becomes $u^{\star}_{0,1} \approx
1585.3$, and the first four (including multiplicities) vibrational
frequencies in Fig.~\ref{fig:NoHoles} have been tuned out. For four
equally spaced clamping points, the lowest mode would have been
$\lambda^{\star}_{2,0} \approx1218.9$ resulting in the lowest three
modes being tuned out.

%%%%%%%%%%%%%%%%%%%%%%%%%%%%%%%%%%%%%%%%%%%%%%%%%%%%%%%%%%%%%%%%%%%%%%%%

\subsection{Rectangular Domain}\label{sec:rectangle}

\begin{figure}[htbp]
\centering
\subfigure[The first two unclamped modes.]{\includegraphics[height =
2in]{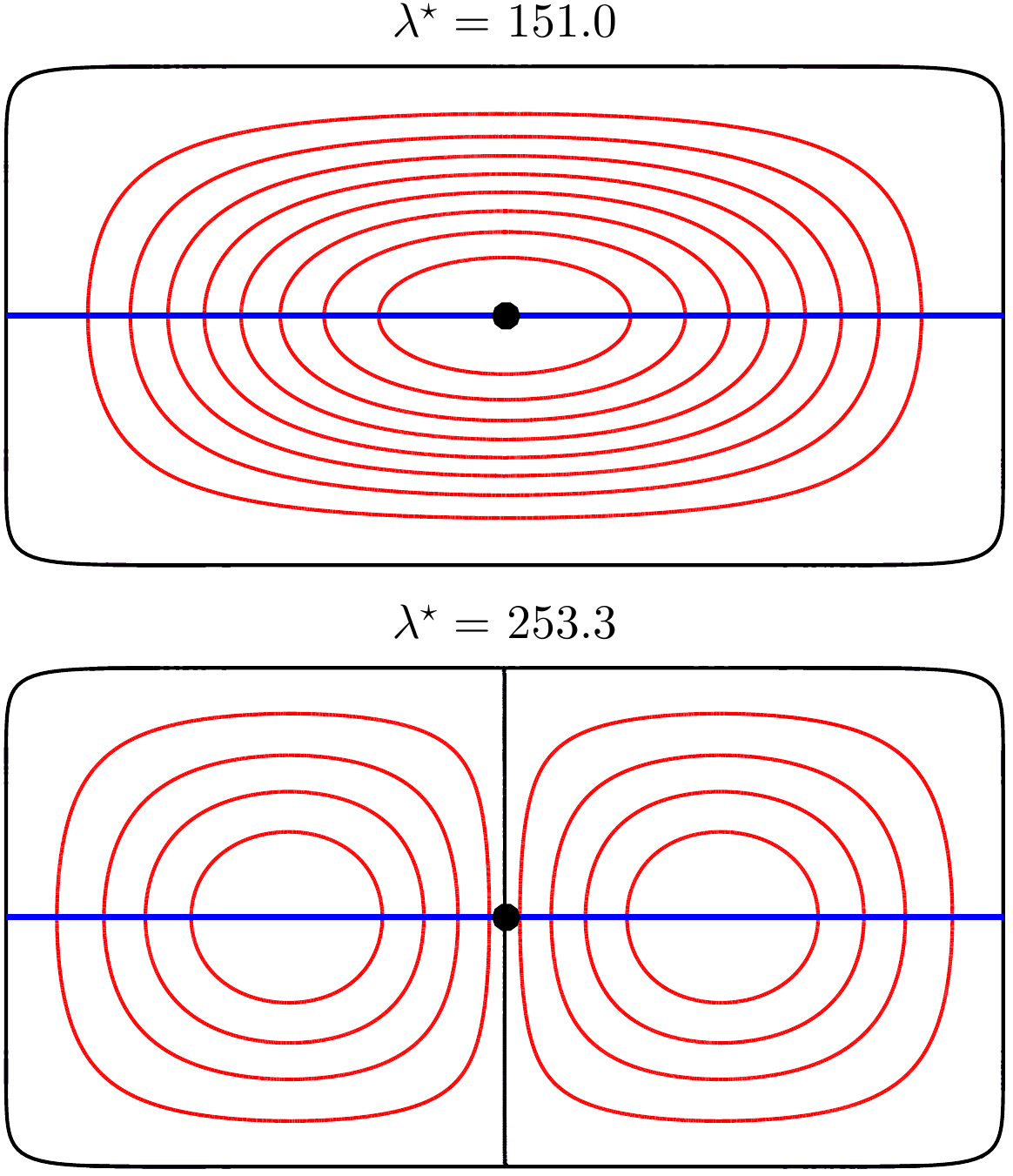} \label{fig:rectangles_a}}\qquad
\subfigure[The first clamped mode along the
centerline]{\includegraphics[height = 2in]{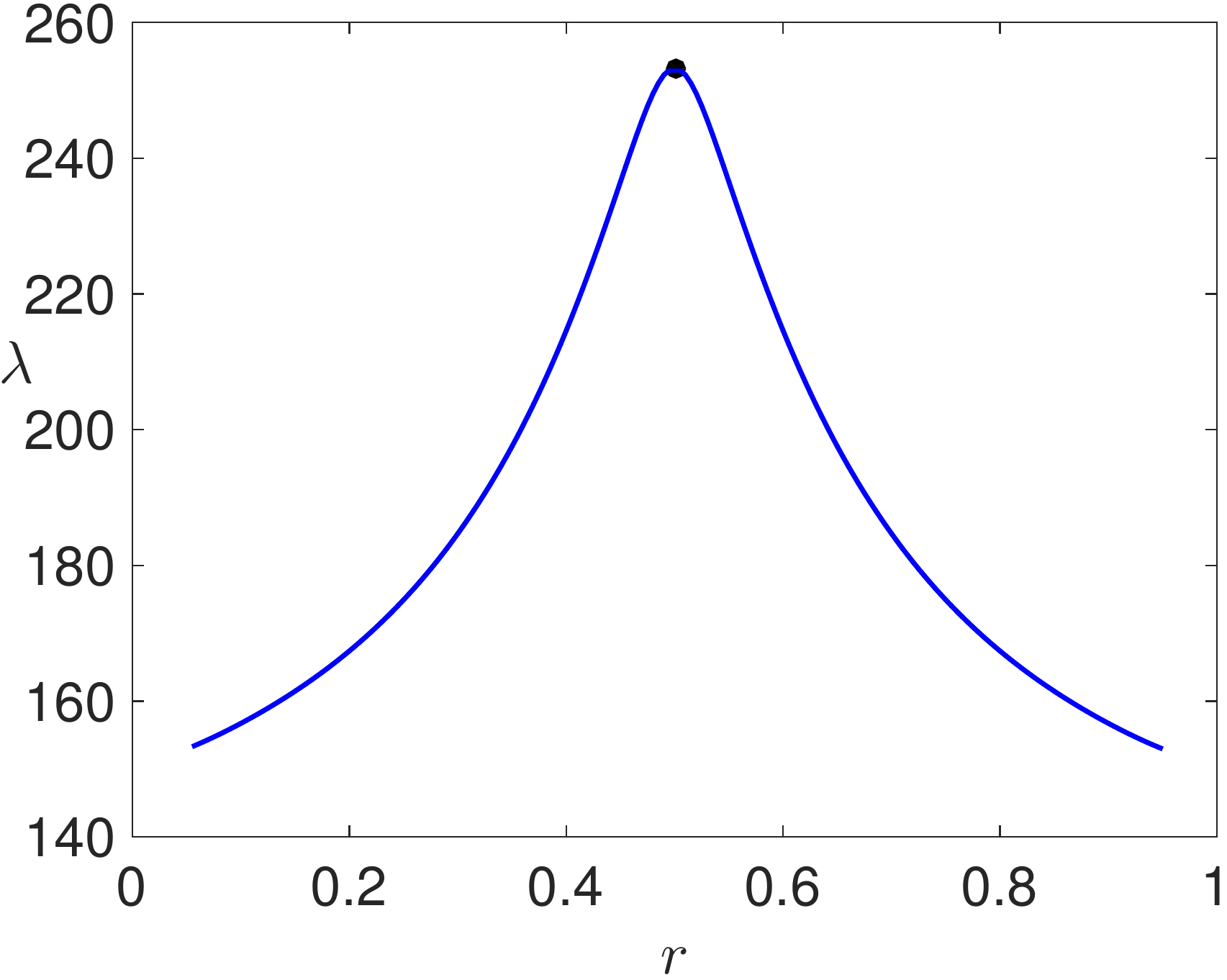} \label{fig:rectangles_b} }
\parbox{0.75\textwidth}{\caption{Results for the rectangular
domain~\eqref{eqn:rectangle} with $a=\sqrt{2}$, $b=1/\sqrt{2}$ and a
single clamped point on the horizontal axis. The eigenvalue attains its
maximum when clamping is at the origin.
\label{fig:rectangles}}}
\end{figure}
In this section, we consider rectangular domains and demonstrate a
qualitative agreement with a previous study~\cite{FM} on the
localization of eigenfunctions of~\eqref{eqn:intro}.  We parameterize
the boundary as
\begin{equation}\label{eqn:rectangle}
  \partial\Omega = (x,y) = (a\,r(\theta) \cos\theta, b\,r(\theta) \sin\theta), \qquad
  \theta \in [0,2\pi), \qquad
  r(\theta) = \left[ \cos^{p}\theta  + \sin^{p} \theta \right]^{-\frac{1}{p}},
\end{equation}
where $p$ is a parameter which regularizes the corners of the rectangle
while $a$ and $b$ give the aspect ratio. We use $p=16$ in our
calculations. Our first simulation investigates the effect of a single
clamping location on the eigenvalue $\lambda$ of~\eqref{eqn:intro}.  We
take a single clamped point and vary its location on the horizontal axis
while calculating $\lambda(r)$ for $r\in(0,1)$ where $r=0$ and $r=1$
correspond to the left and right hand boundaries. In the curve
Fig.~\ref{fig:rectangles_b} we see that the lowest eigenvalue attains a
maximum at $r=1/2$ when the clamping occurs at the origin. The value of
$\lambda$ at that peak corresponds to the second unclamped mode seen in
eigenfunction of Fig.~\ref{fig:rectangles_a}. Therefore a single
clamping point placed at the center of the rectangle completely removes
the first mode from the spectrum of~\eqref{HoleFree}.
\begin{figure}[htbp]
\centering
\subfigure[One clamped point.]{\includegraphics[width= 0.45\textwidth]{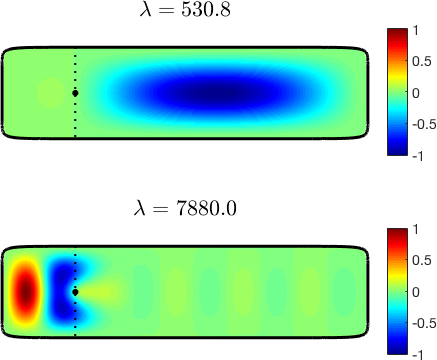} \label{fig:ClampedRectangles_a}}\qquad
\subfigure[Two clamped points.]{\includegraphics[width= 0.45\textwidth]{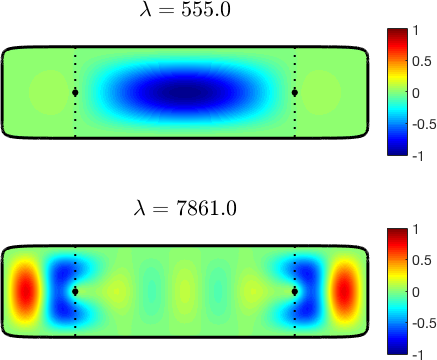} \label{fig:ClampedRectangles_b} }
\parbox{0.75\textwidth}{\caption{Mode confinement effect with one (left
panel) and two (right panel) point constraints on higher (bottom) and
lower (top) modes of~\eqref{eqn:intro} in the rectangular
domain~\eqref{eqn:rectangle} with $a=2$, $b= 1/2$. Clamping occurs on
the horizontal axis at locations $\pm(1.2,0)$.  The simulations use
$N=512$ discretization points of the domain's boundary.
\label{fig:ClampedRectangles} }}
\end{figure}

In our second simulation, we investigate the effect clamping points has
on the eigenfunctions oscillations (cf.~\cite{FM}). In
Fig.~\ref{fig:ClampedRectangles} we display two modes
for~\eqref{eqn:intro} for the rectangular domain with $a=2$, $b = 1/2$
with one and two clamped points. The main observation is that a result
of clamping is a confinement region, i.e.~the eigenfunction is
effectively zero on either the left or right of the clamping point.

\subsection{Elliptical Domain}
\begin{figure}[htbp]
\centering
\includegraphics[width = 0.95\textwidth]{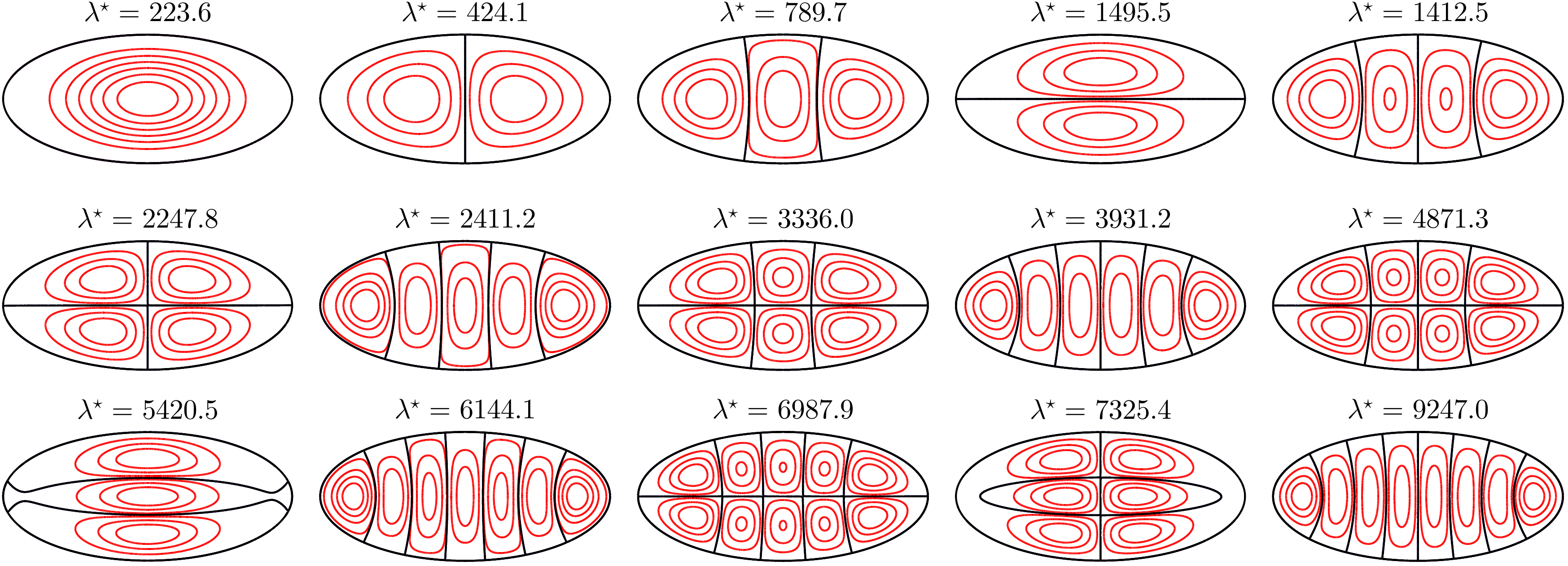}
\parbox{0.75\textwidth}{\caption{The first $15$ clamped modes
of~\eqref{HoleFree} of an ellipse~\eqref{def:ellipse} with $a= 3/2$, $b
= 2/3$. The nodal lines ($u^{\star}=0$) are in black.
\label{fig:NoHolesEllipse}}}
\end{figure}

In this section, we consider equation~\eqref{eqn:intro} on the
elliptical domain defined parametrically as
\begin{equation}\label{def:ellipse}
\partial\Omega = (x,y) = (a \cos\theta,   b \sin\theta), \qquad
  \theta \in [0,2\pi).
\end{equation}
As with the disk problem, the defect free eigenfunctions $u^{\star}$
satisfying equation~\eqref{HoleFree} gives insight into the effect that
puncturing will have on the modes. In Fig.~\ref{fig:NoHolesEllipse}, we
display the first $15$ modes of~\eqref{HoleFree} for $a = \frac{3}{2}$,
$b = \frac{2}{3}$ with nodal lines.  The choice $a = b^{-1}$ is made so
that the ellipse has area $\pi$.

\begin{figure}[htbp]
\centering
\subfigure[Ellipse puncture patterns.]{\includegraphics[width =
0.475\textwidth]{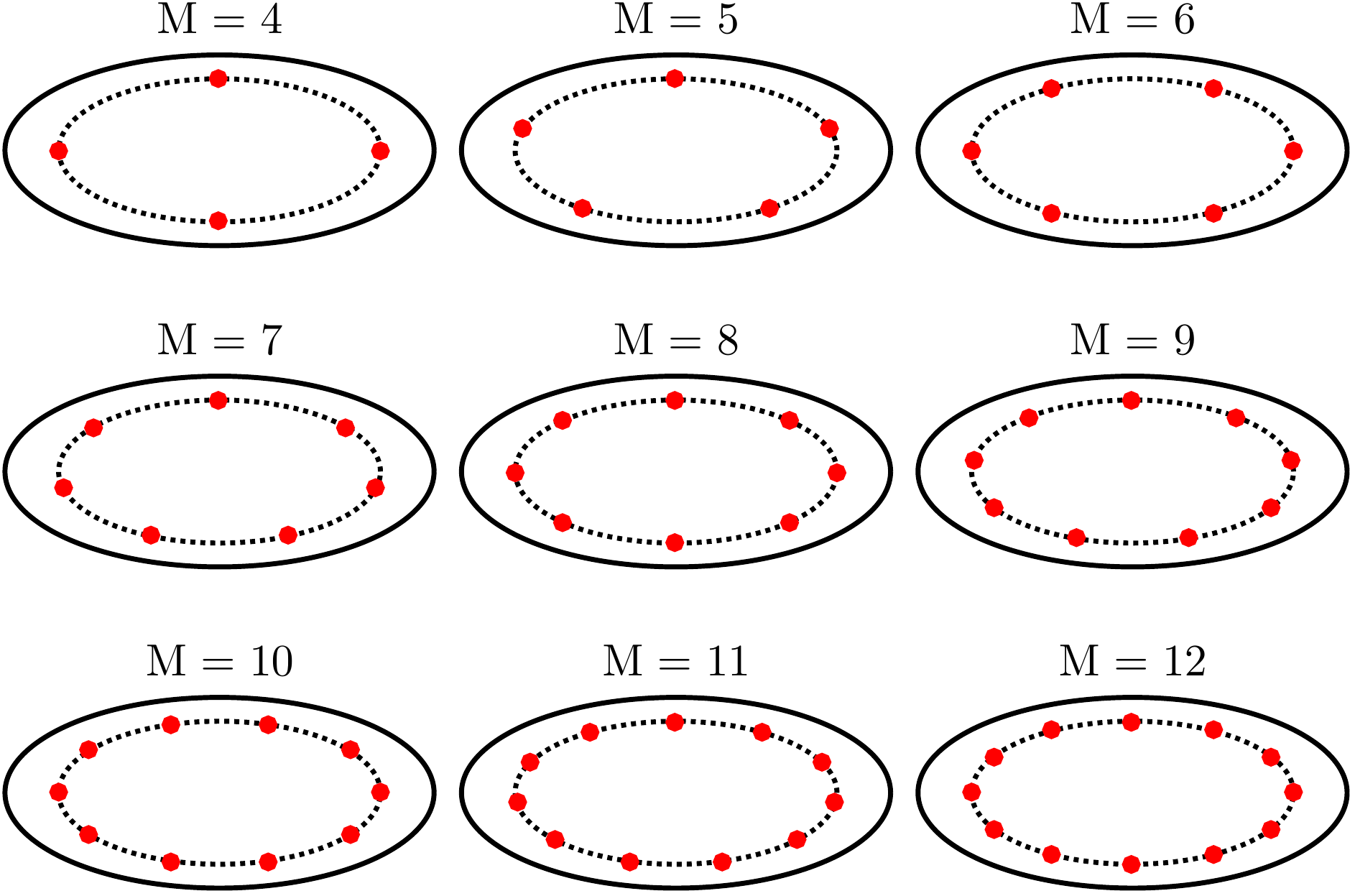}}\qquad
\subfigure[Rectangular puncture patterns.]{\includegraphics[width =
0.475\textwidth]{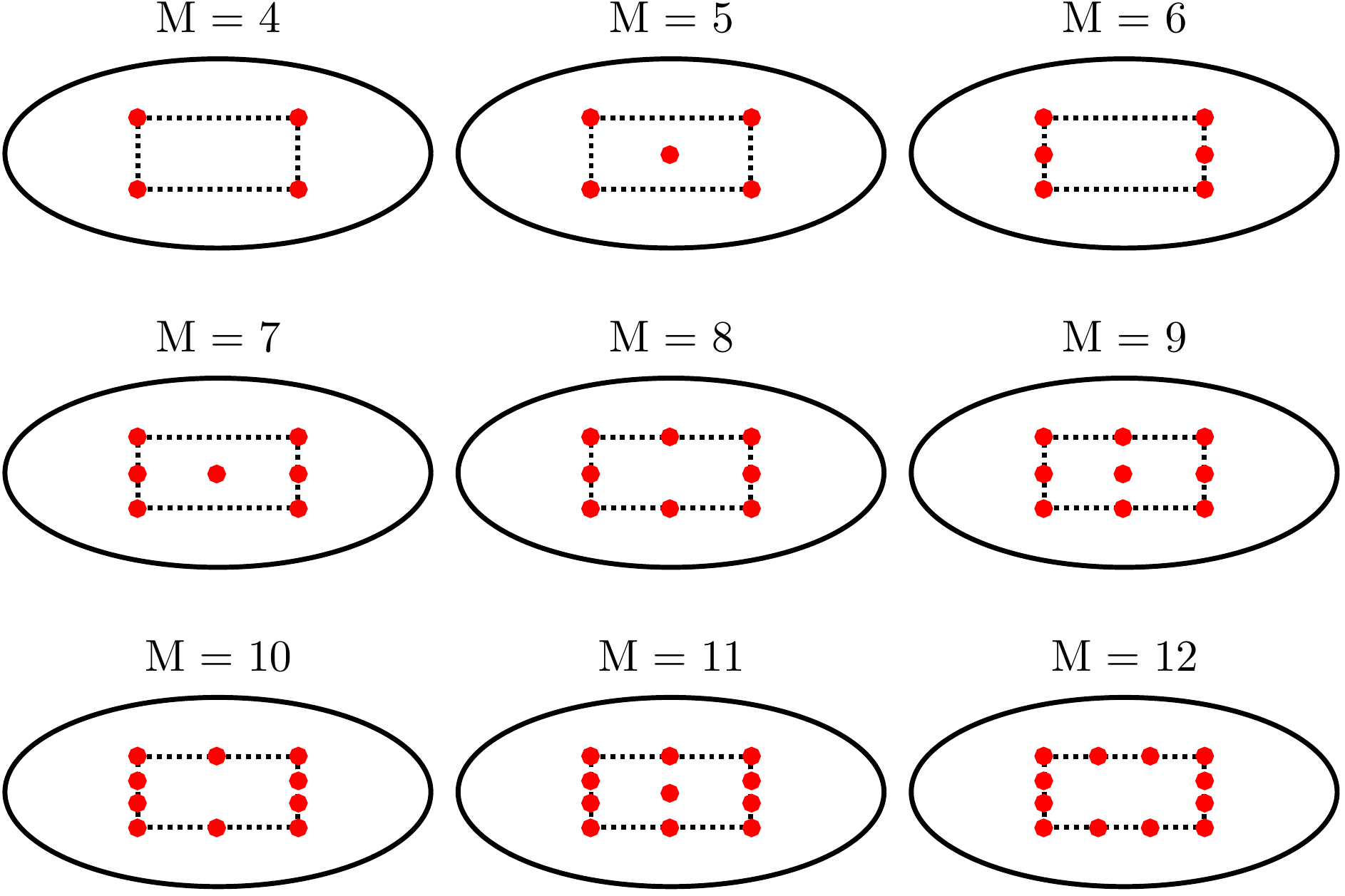}}
\parbox{0.75\textwidth}{\caption{Elliptical and rectangular puncture
patterns for $M=4,\ldots,12,$ inside an ellipse \eqref{def:ellipse} with
$a = 3/2$, $b = 2/3$. The aspect ratio of the patterns matches
that of the outer ellipse. \label{fig:EllipsePatterns} }}
\end{figure}

As with the disk domain, we determine solutions of~\eqref{eqn:intro}
over one parameter families of puncture patterns. In this example we
choose inscribed ellipses, circular patterns, and rectangular patterns.
For the elliptical and rectangular cases, the aspect ratio of the
puncture pattern is chosen to be the same as the outer ellipse and the
pattern is stretched in a uniform way by a single parameter $r$
(cf.~Fig.~\ref{fig:EllipsePatterns}).

\begin{figure}[htbp]
\centering
\includegraphics[width = \textwidth]{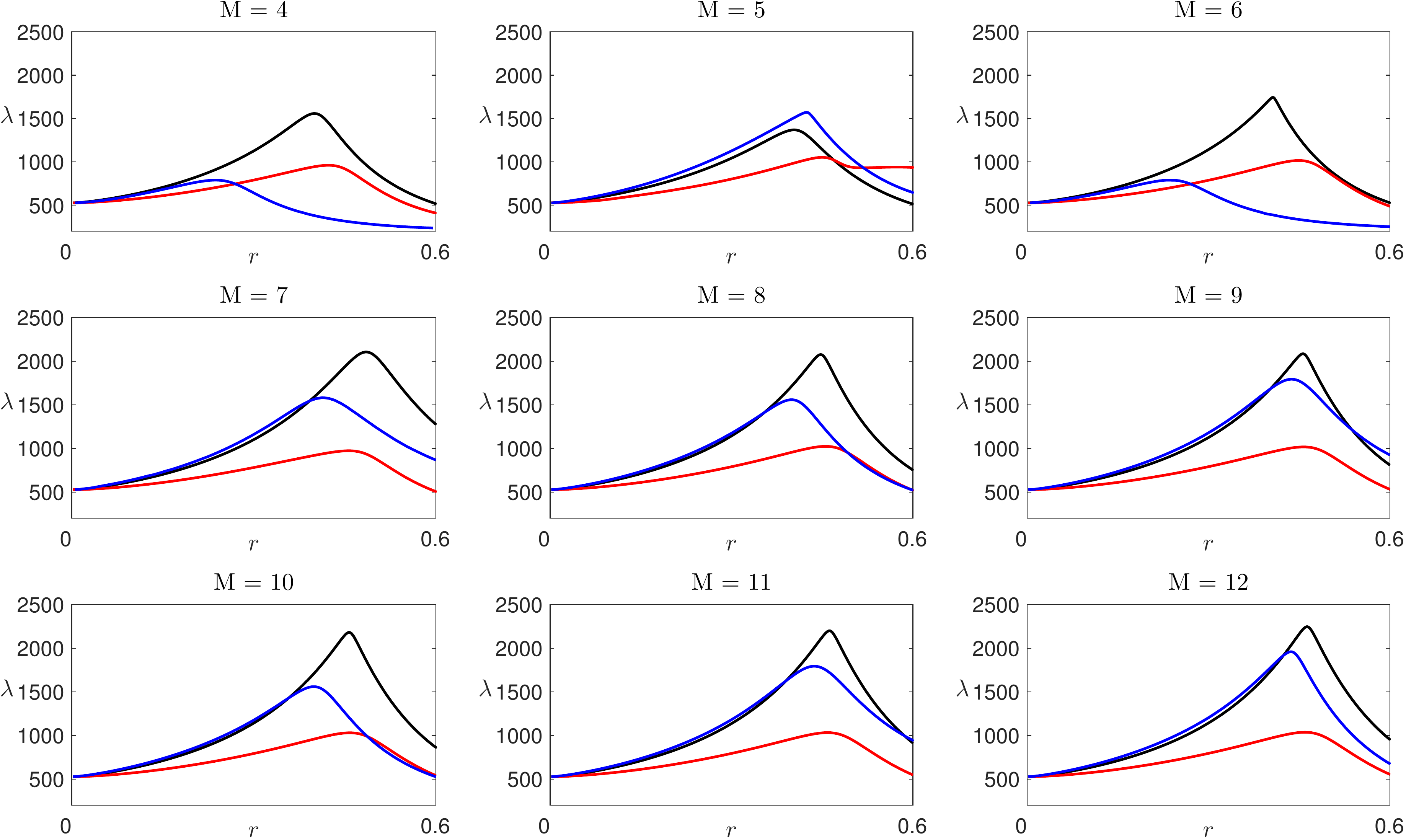}\
\parbox{0.75\textwidth}{\caption{The first eigenvalue of \eqref{eqn:intro} for the ellipse \eqref{def:ellipse} with $a= 3/2$, $b = 2/3$ and $M=4,\ldots,12,$ punctures. Curves for rectangular (blue), elliptical (black) and circular (red) patterns. \label{fig:EllipseResults} }}
\end{figure}

The results for the puncture patterns corresponding to
Fig.~\ref{fig:EllipsePatterns} and the circular pattern
Fig.~\ref{fig:SingleRing} are given in Fig.~\ref{fig:EllipseResults}.
The curves show that the pattern that generates the maximum eigenvalue
varies for different $M$. The horizontal axis is a positive parameter
$r$ which is a scale factor controlling the size of the pattern. The
elliptical puncture pattern generates the largest eigenvalue for all $M$
tested except for $M=5$ where the rectangular patter generates a higher
value. The circular pattern generally results in a lower maximum
eigenvalue than the other patterns tested, expect in the $M=4$ case
where the circular example generates a slightly larger value than the
rectangle.

The discussion from the previous circular and rectangular examples,
together with formula~\eqref{EigDifference}, suggest that clamping
along the nodal lines of the unperturbed mode results in a maximum
deviation of the eigenvalue $\lambda$ of~\eqref{eqn:intro}. This
suggests that the better performance of the elliptical pattern in
maximizing the eigenvalue is that it more closely places punctures on
the nodal lines (cf.~Fig.~\ref{fig:EllipsePatterns}) of
$u^{\star}(\xx)$.

%

%%%%%%%%%%%%%%%%%%%%%%%%%%%%%%%%%%%%%%%%%%%%%%%%%%%%%%%%%%%%%%%%%%%%%%%%
\subsection{Non-symmetric geometry}\label{sec:nonsymmetric}

\begin{figure}[htbp]
\centering
\subfigure[The first and second unpunctured
modes.]{\includegraphics[width=0.5\textwidth]{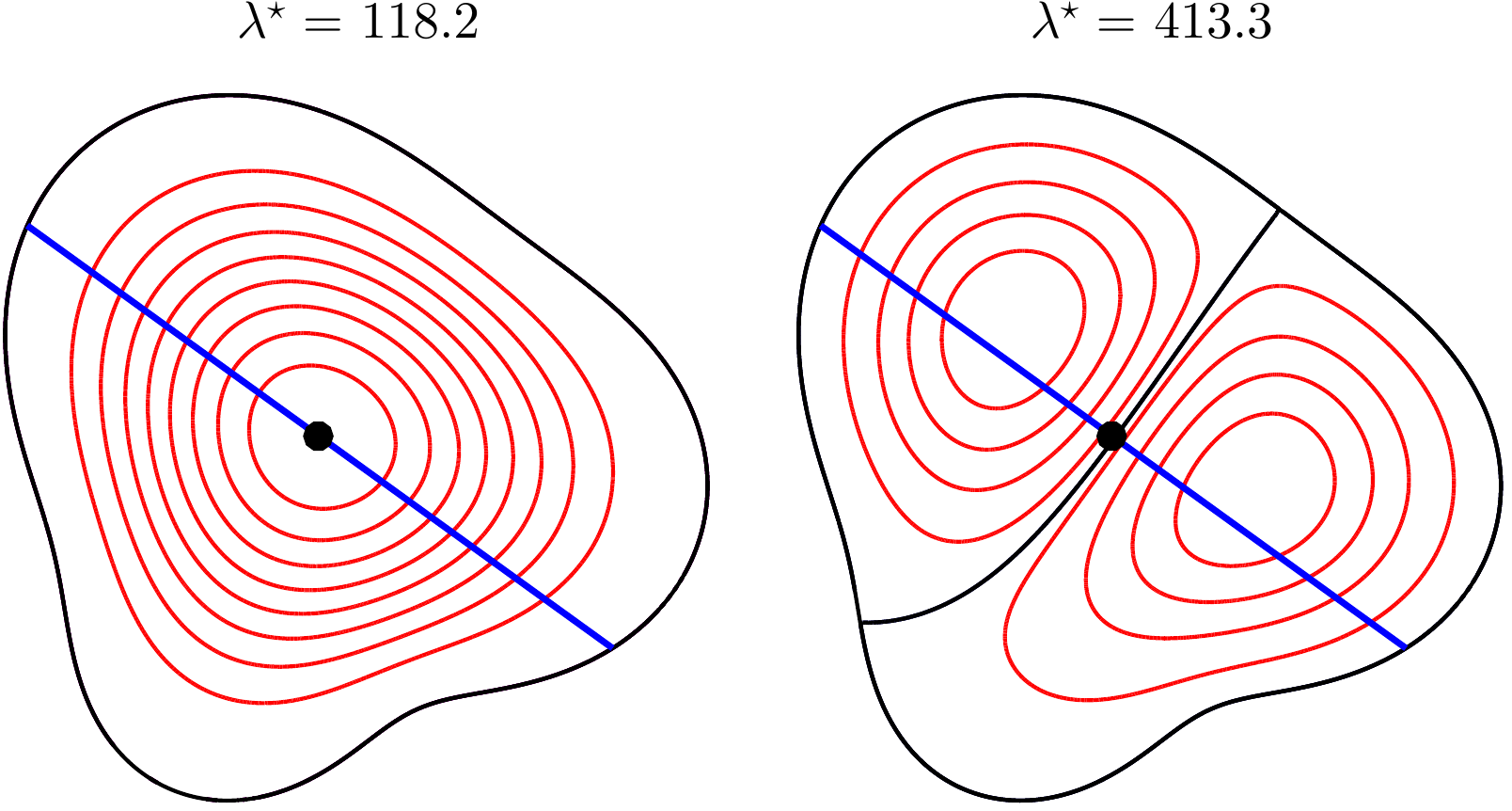}\label{fig:asym_a}}\qquad
\subfigure[The first punctured
eigenvalue.]{\includegraphics[width=0.32\textwidth]{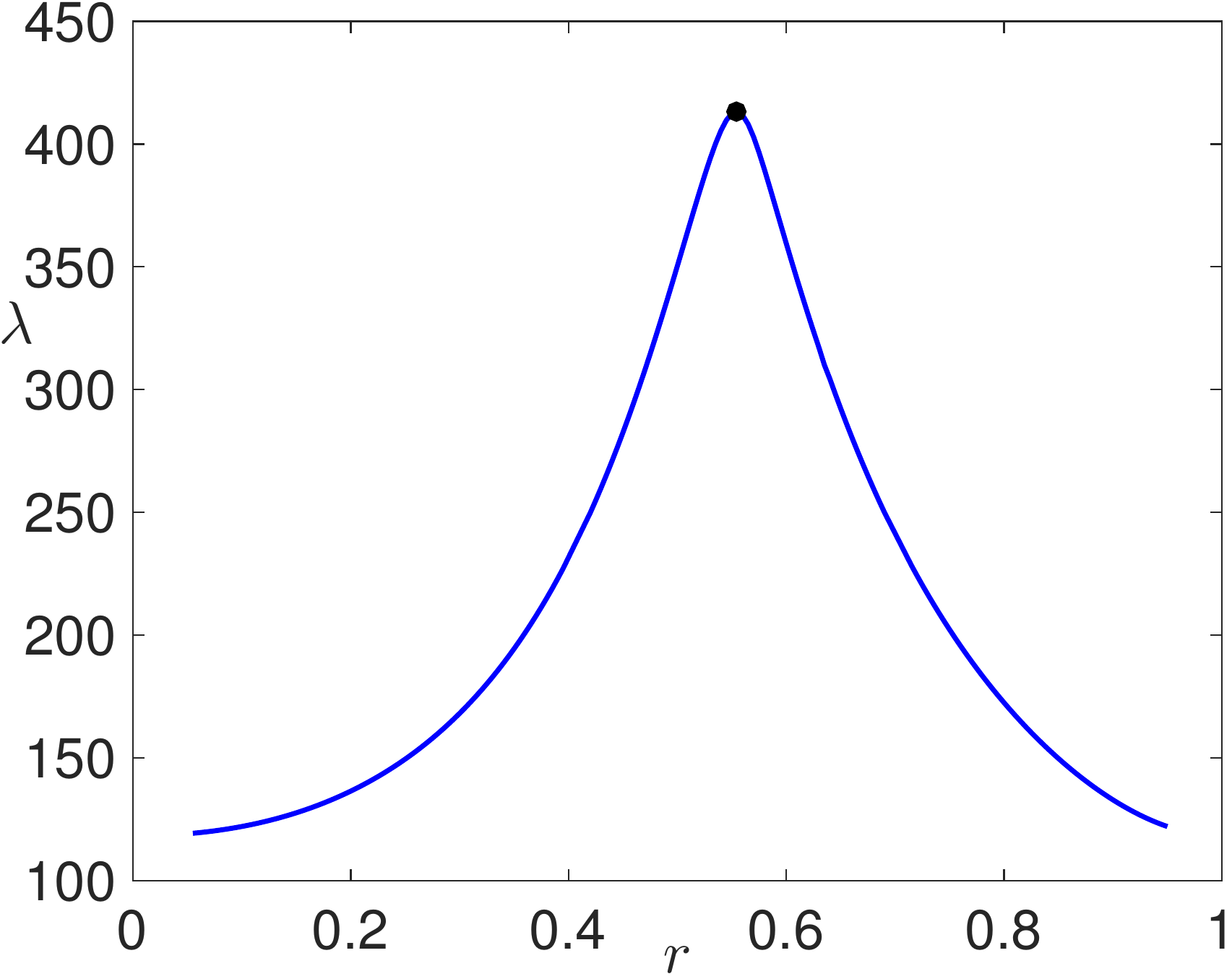}\label{fig:asym_c}}
\parbox{0.75\textwidth}{\caption{Panel (a): The first two modes of the hole free problem~\eqref{HoleFree}. Panel (b): The first eigenvalue of~\eqref{eqn:intro} for a clamped point located a distance $r$ along the blue line in panel (a). \label{fig:asym} }}
\end{figure}

Here we consider the asymmetric domain whose boundary $\partial\Omega$
is specified parametrically as
\begin{equation*}
%  \label{asydomain}
  \partial\Omega = (r(\theta)\cos\theta,r(\theta)\sin\theta), \qquad 
  \theta\in[0,2\pi); \qquad 
  r(\theta) = 1 + 0.25\sin \theta + 0.15\cos 3\theta,
\end{equation*}
and investigate the solution of~\eqref{eqn:intro} with a single clamped
point. As in the previous examples, the first step is to consider the
eigenvalues of the hole free problem~\eqref{HoleFree}. The first two
modes are shown in Fig.~\ref{fig:asym_a}.

To initiate Newton iterations for this example, we begin with a single
clamped point near the boundary so that the unperturbed mode
$\lambda^{\star}\approx 118.3$ in Fig.~\ref{fig:asym_a} provides a good
initial guess for the system. From this start point, we vary the
puncture location along a straight line (blue line in
Fig.~\ref{fig:asym_a}) through the center of the domain to the opposing
boundary. The maximum of this curve is attained when the clamping
location coincides with the maximum of the first eigenfunction,
occurring at the solid dot shown in Fig.~\ref{fig:asym_a}.
In this example we again see that the first frequency can be tuned out of the vibrational characteristics of the plate by careful placement of just one clamping point.
%
%%%%%%%%%%%%%%%%%%%%%%%%%%%%%%%%%%%%%%%%%%%%%%%%%%%%%%%%%%%%%%%%%%%%%%%%
\subsection{Non-simply connected geometry}
\begin{figure}[htbp]
\centering
\subfigure[Domain and first unclamped mode.]{\includegraphics[width = 0.3\textwidth,clip]{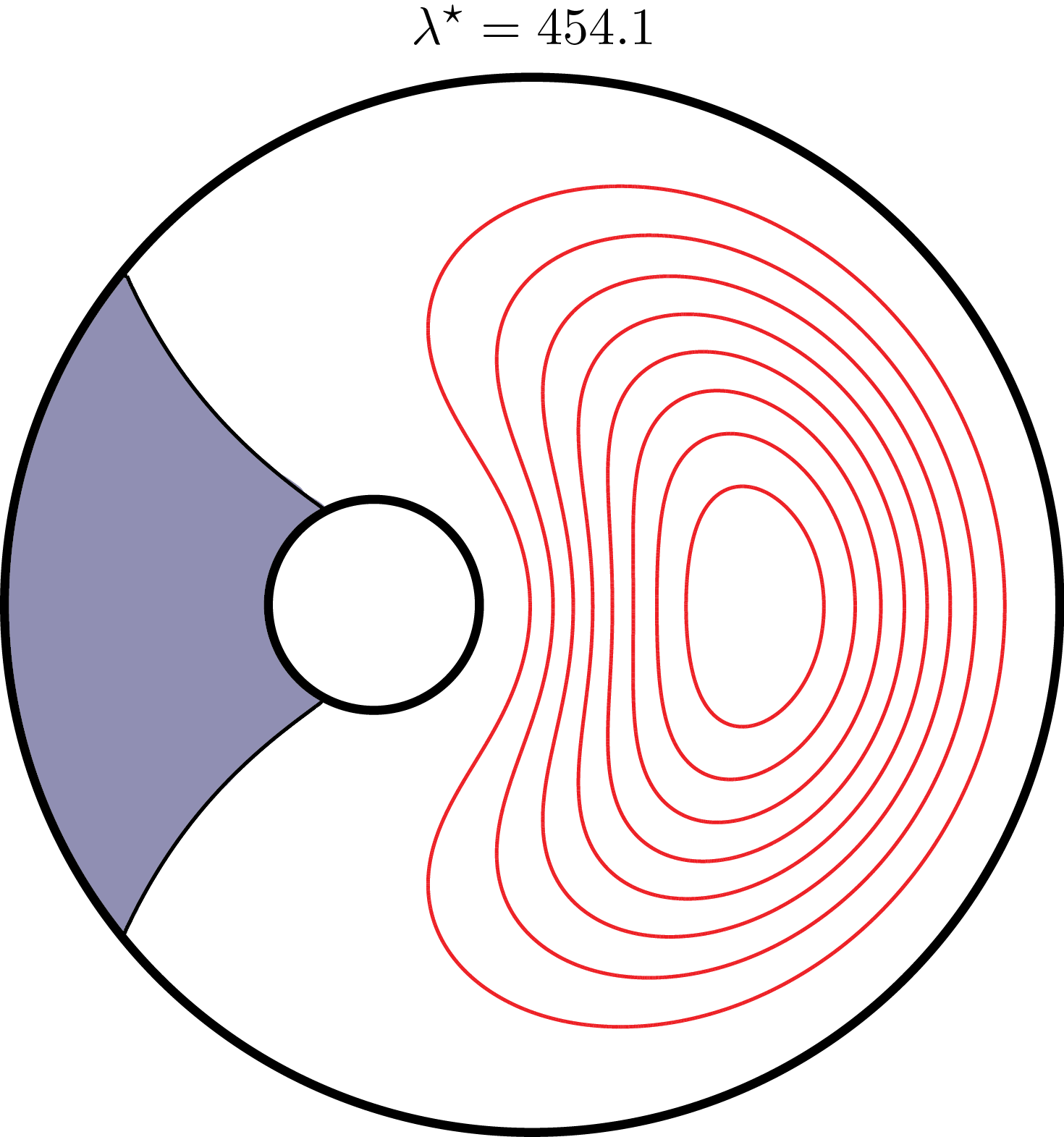} \label{fig:couette1a}}\hspace{1in}
\subfigure[Effect of clamping along an ellipse.]{\includegraphics[width
= 0.4\textwidth,clip]{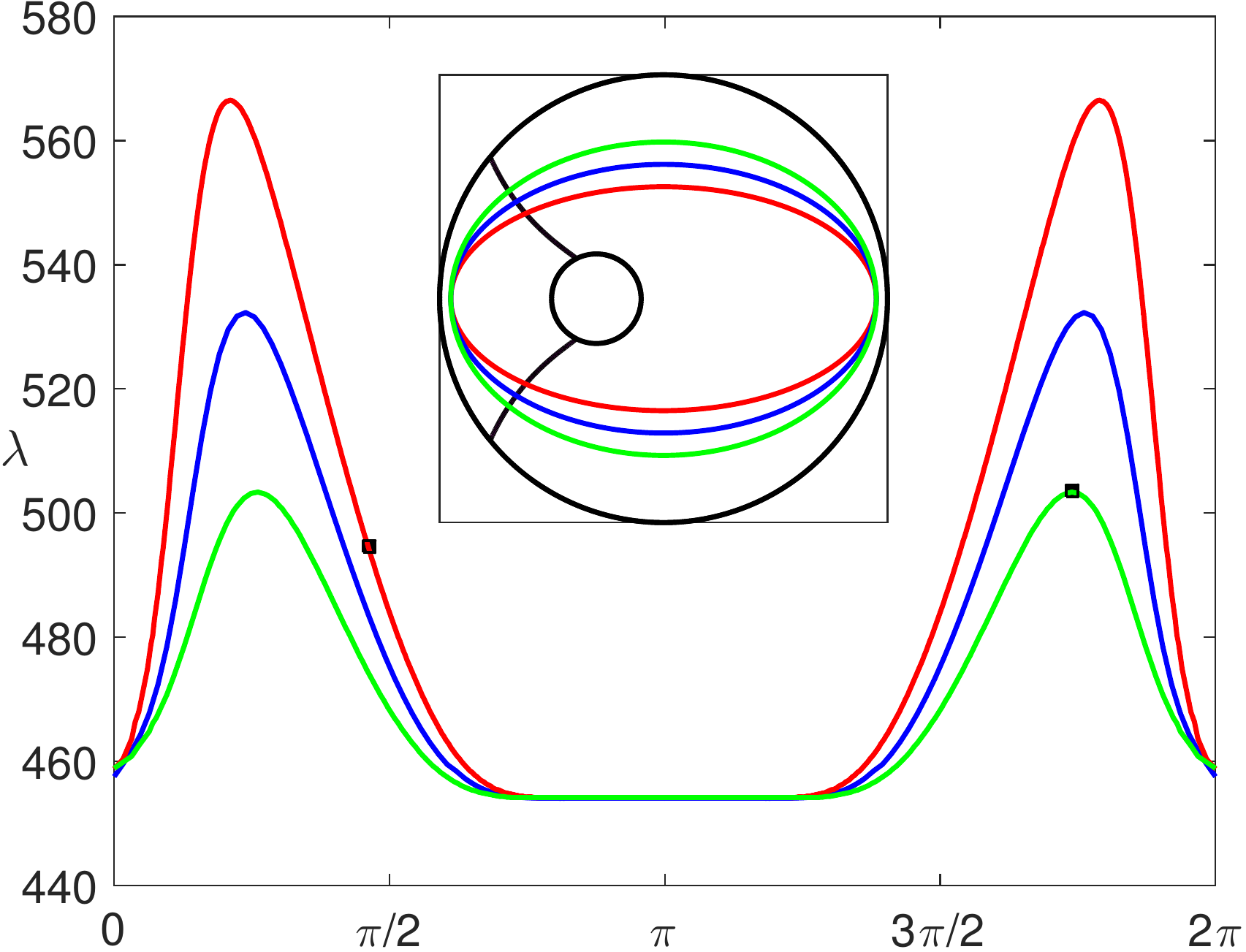} \label{fig:couette1b}}
\parbox{0.75\textwidth}{\caption{Non simply-connected domain example: Left: The
first mode corresponding to $\lambda^{\star} = 454.1$ with the shaded
confinement zone in which
$u^{\star} \approx 0$. Right: The clamped eigenvalue for a single point
on inner ellipses with $a = 0.95$ and $b = 0.5$ (red curve), $0.6$ (blue curve), $0.7$ (green curve). The clamped
modes corresponding to the black squares are shown in Fig.~\ref{fig:couetteModes}.
\label{fig:couette1}}}
\end{figure}

In this section, we consider a non-simply connected domain comprised of
the unit disk with a circle of radius $0.2$ and center $(-0.3,0)$
removed.  This multiply-connected domain has a confinement zone in the
eigenfunction (see Fig.~\ref{fig:couette1}). Throughout this confinement
zone (shaded region of Fig.~\ref{fig:couette1a}), the eigenfunction is
very close to zero indicating again that small defects (holes, clamping
points) cause a global perturbation to the modal characteristics in the
fourth-order problem \eqref{eqn:intro}.

\begin{figure}[htbp]
\centering
\subfigure[$b = 0.5$.]{\includegraphics[width = 0.4\textwidth,clip]{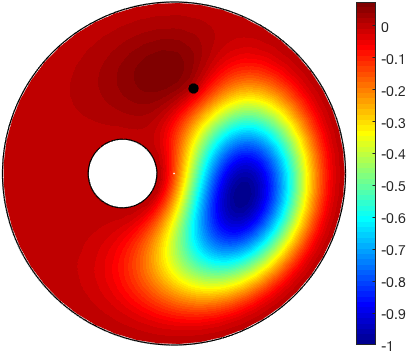} \label{fig:couetteModesa}}\hspace{0.5in}
\subfigure[$b = 0.7$.]{\includegraphics[width = 0.4\textwidth,clip]{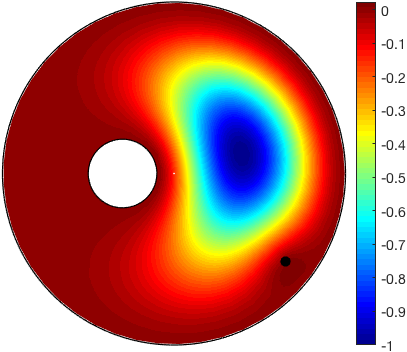} \label{fig:couetteModesb}}
\parbox{0.75\textwidth}{\caption{The clamped modes corresponding to
solid squares in Fig.~\ref{fig:couette1b}. The mode changes sign in the
vicinity of the clamping point. Left: The single clamped point is $\xx_1
= (0.114, 0.496)$ and $\lambda = 494.5$. Right: The single clamped point
is $\xx_1 = (0.648, -0.512)$ and $\lambda = 503.4$.}
\label{fig:couetteModes} }
\end{figure}

To investigate the effect of clamping, we vary the location of a single
clamped point along an inner ellipse with parameterization
\eqref{def:ellipse} for $a=0.95$ and the three values $b= 0.5$, $0.6$,
$0.7$. The results, displayed in Fig.~\ref{fig:couette1b}, show the
eigenvalue as a function of angular position on the ellipse. The flat
region in the center of each curve corresponds to the transit of the
clamped point through the confinement zone. In this region the
eigenfunction is very close to zero ($u^{\star}\approx 0$), therefore
the formula~\eqref{EigDifference} shows that the eigenvalue should
correspond to the unperturbed mode of Fig.~\ref{fig:couette1a}.

In Fig.~\ref{fig:couetteModes} we show two typical clamped modes which correspond to the particular clamping locations marked in Fig.~\ref{fig:couette1b}. In each case a large confinement zone is seen to the left of the hole while the mode is also suppressed in the vicinity of the clamped point.

%%%%%%%%%%%%%%%%%%%%%%%%%%%%%%%%%%%%%%%%%%%%%%%%%%%%%%%%%%%%%%%%%%%%%%%%
\section{Conclusions}\label{sec:conclusions}

This paper has analyzed the modes of vibrations of thin elastic plates
with multiple point constraints. We have developed and validated a novel
boundary integral method for determining the eigenvalues of the
fourth-order bi-Laplacian problem~\eqref{eqn:intro} with multiple
clamped point constraints. This method is third-order accurate and can
be easily applied to a variety of symmetric, asymmetric and
multiply-connected geometries in two dimensions.

Our results indicate that the number and location of clamping points has
two profound effects on the modes of vibrations of thin plates. First,
by placing the clamping locations at the nodal lines of the unclamped
eigenfunctions, certain eigenvalues can be removed from the spectrum of
the problem.  This implies that the vibrational characteristics of the
plate can be manipulated or tuned by judicious placement of a small
number clamping sites. A particularly important consequence of this
effect in engineering applications is that undesirable frequencies of
vibration can be completely removed.  Second, clamping can have the
effect of partitioning the plate into multiple subdomains in which
vibrational modes are largely confined to a subset of those smaller
spatial regions.  This localization effect, previously seen only in
rectangular plates with a single clamping location~\cite{FM}, has been
observed here in other geometries and for multiple clamping locations. 

There are many possibilities for future work that arise from this study.
The order of the numerical method could be improved by adopting
quadrature methods (cf.~\cite{alp1999, kap-rok1997}) suited to integrals
with logarithmic singularities.  In addition, the boundary integral
equation could be more efficiently solved with a fast summation method
such as the kernel-independent fast multipole
method~\cite{yin-bir-zor2004}.  A more significant challenge is to
develop a numerical method for the point constraint eigenvalue
problem~\eqref{eqn:intro} which does not require the solution of a
nonlinear system, such as~\eqref{eqn:mismatch}, but still enforces the
weak logarithmic singularity on the eigenfunction. This would eliminate
the use of Newton's method and the need to carefully control the initial
conditions of the system. Such a method would easily accommodate a much
larger number of clamping locations and allow for more reliable
evaluation of the high frequency modes of~\eqref{eqn:intro}.

%%%%%%%%%%%%%%%%%%%%%%%%%%%%%%%%%%%%%%%%%%%%%%%%%%%%%%%%%%%%%%%%%%%%%%%%

\section*{Acknowledgements} A.E.L.~was supported by the National Science
Foundation under grant DMS 1516753.  B.Q.~was supported by startup funds
at Florida State University.

\begin{appendices}
%%%%%%%%%%%%%%%%%%%%%%%%%%%%%%%%%%%%%%%%%%%%%%%%%%%%%%%%%%%%%%%%%%%%%%%%
\section{Kernels}
\label{app:kernels}
%%%%%%%%%%%%%%%%%%%%%%%%%%%%%%%%%%%%%%%%%%%%%%%%%%%%%%%%%%%%%%%%%%%%%%%%

The four kernels that appear in the integral equation~\eqref{eqn:bie}
are $G_{11} = G_{1}$, $G_{12} = G_{2}$, $G_{21} =
\partial{\nn_{\xx}}G_{1}$, and $G_{22} = \partial{\nn_{\xx}}G_{2}$.
Expressions for $G_{11}$ and $G_{12}$ are in Section~\ref{sec:kernels}.
Here we compute their normal derivatives with respect to the target
point.

\begin{align*}
  G_{21} = -\frac{1}{4\pi\mu^{2}}&\left(
    -3\mu^{4}K_{0}(\mu\rho)\frac{(\rdotn)(\rdotnx)}{\rho^{2}}
    -12\mu^{3}K_{1}(\mu\rho)\frac{(\rdotn)(\rdotnx)}{\rho^{3}}
    +3\mu^{3}K_{1}(\mu\rho)\frac{(\ndotnx)}{\rho}
  \right. \\
  &\left.
    +2\mu^{4}K_{0}(\mu\rho)\frac{(\rdotn)^{3}(\rdotnx)}{\rho^{4}}
    +16\mu^{3}K_{1}(\mu\rho)\frac{(\rdotn)^{3}(\rdotnx)}{\rho^{5}}
    -6\mu^{3}K_{1}(\mu\rho)\frac{(\rdotn)^{2}(\ndotnx)}{\rho^{3}}
  \right. \\
  &\left.
    +6\mu^{2}K_{0}(\mu\rho)\frac{(\ndotnx)}{\rho^{2}} 
    -24\mu^{2}K_{0}(\mu\rho)\frac{(\rdotn)(\rdotnx)}{\rho^4}
    -24\mu^{2}K_{0}(\mu\rho)\frac{(\rdotn)^{2}(\ndotnx)}{\rho^{4}}
  \right. \\
  &\left.
    +48\mu^{2}K_{0}(\mu\rho)\frac{(\rdotn)^{3}(\rdotnx)}{\rho^{6}}
    +96\mu K_{1}(\mu\rho)\frac{(\rdotn)^{3}(\rdotnx)}{\rho^{7}}
    -48\mu K_{1}(\mu\rho)\frac{(\rdotn)^{2}(\ndotnx)}{\rho^{5}}
  \right. \\
  &\left.
     -48\mu K_{1}(\mu\rho)\frac{(\rdotn)(\rdotnx)}{\rho^{5}}
     +12\mu K_{1}(\mu\rho)\frac{(\ndotnx)}{\rho^{3}}
  \right) \\
  -\frac{1}{8\mu^{2}}&\left(
    -3\mu^{4}Y_{0}(\mu\rho)\frac{(\rdotn)(\rdotnx)}{\rho^{2}}
    +12\mu^{3}Y_{1}(\mu\rho)\frac{(\rdotn)(\rdotnx)}{\rho^{3}}
    -3\mu^{3}Y_{1}(\mu\rho)\frac{(\ndotnx)}{\rho}
  \right. \\
  &\left.
    +2\mu^{4}Y_{0}(\mu\rho)\frac{(\rdotn)^{3}(\rdotnx)}{\rho^{4}}
    -16\mu^{3}Y_{1}(\mu\rho)\frac{(\rdotn)^{3}(\rdotnx)}{\rho^{5}}
    +6\mu^{3}Y_{1}(\mu\rho)\frac{(\rdotn)^{2}(\ndotnx)}{\rho^{3}}
  \right. \\
  &\left.
    -6\mu^{2}Y_{0}(\mu\rho)\frac{(\ndotnx)}{\rho^{2}} 
    +24\mu^{2}Y_{0}(\mu\rho)\frac{(\rdotn)(\rdotnx)}{\rho^4}
    +24\mu^{2}Y_{0}(\mu\rho)\frac{(\rdotn)^{2}(\ndotnx)}{\rho^{4}}
  \right. \\
  &\left.
    -48\mu^{2}Y_{0}(\mu\rho)\frac{(\rdotn)^{3}(\rdotnx)}{\rho^{6}}
    +96\mu Y_{1}(\mu\rho)\frac{(\rdotn)^{3}(\rdotnx)}{\rho^{7}}
    -48\mu Y_{1}(\mu\rho)\frac{(\rdotn)^{2}(\ndotnx)}{\rho^{5}}
  \right. \\
  &\left.
    -48\mu Y_{1}(\mu\rho)\frac{(\rdotn)(\rdotnx)}{\rho^{5}}
    +12\mu Y_{1}(\mu\rho)\frac{(\ndotnx)}{\rho^{3}}
  \right).
\end{align*}
\begin{align*}
  G_{22} = -\frac{1}{4\pi}&\left( 
    8 K_{0}(\mu\rho)\frac{(\rdotn)^2(\rdotnx)}{\rho^{4}}
    +\frac{16}{\mu}K_{1}(\mu\rho)\frac{(\rdotn)^{2}(\rdotnx)}{\rho^{5}}
    -4K_{0}(\mu\rho)\frac{(\rdotn)(\ndotnx)}{\rho^{2}}
    \right. \\
    &\left.
    -\frac{8}{\mu}K_{1}(\mu\rho)\frac{(\rdotn)(\ndotnx)}{\rho^{3}}
    -\mu K_{1}(\mu\rho)\frac{(\rdotnx)}{\rho} 
    -\frac{4}{\mu}K_{1}(\mu\rho)\frac{(\rdotnx)}{\rho^{3}}
    \right. \\
    &\left.
    -2K_{0}(\mu\rho)\frac{(\rdotnx)}{\rho^{2}} 
    +2\mu K_{1}(\mu\rho)\frac{(\rdotn)^{2}(\rdotnx)}{\rho^{3}}
  \right) \\
  +\frac{1}{8}&\left( 
    8 Y_{0}(\mu\rho)\frac{(\rdotn)^2(\rdotnx)}{\rho^{4}}
    -\frac{16}{\mu}Y_{1}(\mu\rho)\frac{(\rdotn)^{2}(\rdotnx)}{\rho^{5}}
    -4Y_{0}(\mu\rho)\frac{(\rdotn)(\ndotnx)}{\rho^{2}}
    \right. \\
    &\left.
    +\frac{8}{\mu}Y_{1}(\mu\rho)\frac{(\rdotn)(\ndotnx)}{\rho^{3}}
    -\mu Y_{1}(\mu\rho)\frac{(\rdotnx)}{\rho} 
    +\frac{4}{\mu}Y_{1}(\mu\rho)\frac{(\rdotnx)}{\rho^{3}}
    \right. \\
    &\left.
    -2Y_{0}(\mu\rho)\frac{(\rdotnx)}{\rho^{2}} 
    +2\mu Y_{1}(\mu\rho)\frac{(\rdotn)^{2}(\rdotnx)}{\rho^{3}}
  \right).
\end{align*}
\end{appendices}

%\bibliography{holesbib}
%\bibliographystyle{plain}
%\bibliographystyle{plainnat}
%\bibliography{refs}
%\biboptions{sort&compress}

\end{document}